\documentclass[preprint,12pt,sort&compress]{elsarticle}
\usepackage{epsfig}
\usepackage{svg}
\usepackage{amssymb,amsmath}
\usepackage{latexsym}
\usepackage{graphicx}
\usepackage{longtable}
\usepackage{amsfonts}
\usepackage{pstricks}
\usepackage{pst-node}
\usepackage{hyperref}
\usepackage{upquote}
\usepackage{arydshln}
\usepackage{mathtools}
\usepackage{algorithm}% http://ctan.org/pkg/algorithm
\usepackage{algpseudocode}%
\usepackage[skip=5pt]{caption}
\usepackage{subcaption}
\usepackage{tikz}
\usepackage{pgfplots}
\usepackage[colorinlistoftodos]{todonotes}
\usepackage{soul}
\usetikzlibrary{shapes.geometric, arrows}
\tikzstyle{startstop} = [rectangle, rounded corners, text centered, draw=black]
\tikzstyle{process} = [rectangle, text centered, draw=black]
\tikzstyle{decision} = [diamond, text centered, draw=black]
\tikzstyle{arrow} = [thick,->,>=stealth]
\hypersetup{
   unicode=false,          % non-Latin characters in Acrobatï¿œs bookmarks
   pdftoolbar=true,        % show Acrobatï¿œs toolbar?
   pdfmenubar=true,        % show Acrobatï¿œs menu?
   pdffitwindow=false,     % window fit to page when opened
   pdfstartview={FitH},    % fits the width of the page to the window
   pdftitle={My title},    % title
   pdfauthor={Author},     % author
   pdfsubject={Subject},   % subject of the document
   pdfcreator={Creator},   % creator of the document
   pdfproducer={Producer}, % producer of the document
   pdfkeywords={keywords}, % list of keywords
   pdfnewwindow=true,      % links in new window
   colorlinks=true,                    % false: boxed links; true: colored links
   linkcolor=blue,                      % color of internal links
   citecolor=red,          % color of links to bibliography
   filecolor=cyan,                     % color of file links
   urlcolor=blue                       % color of external links
}
\usepackage{txfonts}
\usepackage{bm}
\usepackage{setspace}
\usepackage{url}
\usepackage{multicol}
\usepackage[top=1in,bottom=1in,left=1in,right=1in]{geometry}
\graphicspath{{figures/}}
\allowdisplaybreaks
\usepackage{titlesec}

\titleformat{\paragraph}
{\normalfont\normalsize\itshape}{\theparagraph}{1em}{}
\titlespacing*{\paragraph}
{0pt}{3.25ex plus 1ex minus .2ex}{1.5ex plus .2ex}

\newcommand{\framework}[1]{GraPTop}
\let\oldReturn\Return
\renewcommand{\Return}{\State\oldReturn}

\def\d2divergence{\pmb{\nabla}^{\ast 2}}

\def\u_j{{u_j}}

\def\x_j{{x_j}}

\begin{document}
%
%---------------- FRONT MATTER -------------------------------------------
%
\begin{frontmatter}
\title{A geometrically informed algebraic multigrid preconditioned iterative approach for solving high-order finite element systems}

\author[utah1,shu]{Songzhe Xu}
\cortext[cor1]{Corresponding author}
\author[utah2]{Majid Rasouli}
\author[utah1,utah2]{Robert M. Kirby}
\author[eng]{David Moxey}
\author[utah2]{Hari Sundar\corref{cor1}}
\ead{hari@cs.utah.edu} 

\address[utah1]{Scientific Computing and Imaging Institute, University of Utah, USA}
\address[utah2]{School of Computing, University of Utah, USA}
\address[eng]{College of Engineering, Mathematical and Physical Sciences, University of Exeter, UK}
\address[shu]{State Key Laboratory of Advanced Special Steels, School of Materials Science and Engineering, Shanghai University, Shanghai, China}

%\journal{Journal of Computational and Applied Mathematics}
\journal{arXiv}

\begin{abstract}
Algebraic multigrid (AMG) is conventionally applied in a black-box fashion, agnostic to the underlying geometry. In this work, we propose that using geometric information -- when available -- to assist with setting up the AMG hierarchy is beneficial, especially for solving linear systems resulting from high-order finite element discretizations. High-order problems draw considerable interest to both the scientific and engineering communities, but lack efficient solvers, at least open-source codes, tailored for unstructured high-order discretizations targeting large-scale,
real-world applications. %However, to the best of the authors' knowledge, there are no AMG solvers that are tailored for high-order problems, or at least, no such solvers publicly available.
For geometric multigrid, it is known that using $p$-coarsening before $h$-coarsening can provide better scalability, but setting up $p$-coarsening is non-trivial in AMG. We develop a geometrically informed algebraic multigrid (GIAMG) method, as well as an associated high-performance computing program, which is able to set up a grid hierarchy that includes $p$-coarsening at the top grids with minimal information of the geometry from the user. A major advantage of using $p$-coarsening with AMG -- beyond the benefits known in the context of geometric multigrid (GMG) -- is the increased sparsification of coarse grid operators. We extensively evaluate GIAMG by testing on the 3D Helmholtz and incompressible flow problems, and demonstrate mesh-independent convergence, and excellent parallel scalability. We also compare the performance of GIAMG with existing AMG packages, including Hypre and ML.
      
\end{abstract}

\begin{keyword}
high-order or spectral/$hp$ element\sep
geometrically informed algebraic multigrid\sep
$p$-coarsening\sep
preconditioning\sep
high-performance computing
%iterative solver
\end{keyword}

\end{frontmatter}

\tableofcontents

\section{Introduction}
% why we solve for higher order systems
%\par $\boldsymbol{Problem \, statement}$: 
The development of robust, efficient solvers that utilize high-order finite element methods (HO-FEM),
also sometimes classified as spectral/$hp$ element methods, is an area of considerable interest to both the scientific and engineering communities \cite{kopriva2009,karniadakis2013spectral}. 
The use of higher order polynomial expansions within 
elements carries a number of benefits, as seen 
from two main perspectives. Numerically, these methods exhibit far lower levels of numerical dispersion and dissipation at
higher polynomial orders, which makes them a particularly well-suited approximation choice in areas such as computational fluid dynamics (CFD). 
CFD generally solves governing equations involving Navier--Stokes equations, or their variations, to study the fluid-related physics. 
%The accurate time-advection of energetic structures such as vortices is a key concern (e.g., \cite{vincent2014,lombard-2016a}) in CFD, requiring accurate treatment for numerical dispersion and dissipation. 
%
Among the traditional methods for numerically solving fluid mechanics, finite volume method is the most popular approach due to its better performance for the conservation laws, while finite element method rapidly develops in recent decades as approaches (such as Petrov-Galerkin~\cite{Brooks1982}, discontinuous Galerkin~\cite{Bassi1997}, velocity correction~\cite{Karniadakis1991}, etc) have been proposed to successfully address the non-symmetric convection term.
In CFD, the trend of increasing computational power has led to interest in large-eddy simulation (LES) or direct numerical simulation (DNS) at increasing levels of fidelity~\cite{vincent2014}. In this context, the accurate advection of energetic vortices and other fluid structures over large time- and length-scales is a key concern, for example in the simulation of jetting vortices~\cite{lombard-2016a}, and these numerical properties help to preserve and accurately represent such structures on relatively coarse grids that would be otherwise very diffusive at lower orders. This enables a more effective balance to be achieved between accuracy and computational efficiency, whilst the elemental decomposition of the method permits complex geometries such as whole car simulations to be considered~\cite{mengaldo-2021}.
Furthermore, given current hardware trends, another appealing property of these methods has been their computational performance.
Although the cost per degree of freedom in terms of algorithmic floating point operations (FLOPS) increases substantially with polynomial degree, 
the use of higher order expansions leads to formulations of the underlying equations of state that involve dense, 
compact kernels for key finite element operators, such as inner products and derivatives. This is important from the perspective of modern
hardware, where increasingly the bottleneck in performance is memory bandwidth as opposed to the clock speed of processors. 
The underlying arithmetic intensity of the algorithm at hand (i.e., the number of floating-point operations performed
for each memory operation) is therefore key to attaining optimal
performance. This is where high-order methods have a significant advantage over lower order methods (e.g., \cite{fehn-2018,Moxey2019}).

The main challenges in high-order discretizations are that matrices are more dense compared with low-order methods, and that they lose structural properties such as the M--matrix property, which often allows one to prove convergence of iterative solvers.
%\todo{\tiny challenges that we aim to address}
For our target applications, multigrid is a preferred method for solving the resulting symmetric linear systems with optimal and mesh-independent convergence \cite{trottenberg2000multigrid}.
However, when used with unstructured meshes, the use of geometric multigrid (GMG) is not straightforward, requiring the use of algebraic multigrid (AMG) instead. Unlike GMG, where the coarse-grid operators are typically obtained via geometric coarsening and rediscretizations, AMG does not use the geometric information but instead builds the coarse grid operators algebraically. Algebraic coarsening results in a loss of sparsity, leading to poor performance and scalability, especially for large-scale problems. 
This is usually addressed by employing sparsity control techniques, such as dropping values below a specified threshold at coarser levels. Such measures can address the sparsity issues, albeit with some loss of convergence rates. Additionally, they introduce an additional parameter that needs to be tuned, making the overall solver or preconditioner less robust.
This issue is exacerbated for high-order operators that are denser to begin with, and get even denser on coarser grids, leading to poor scalability. 

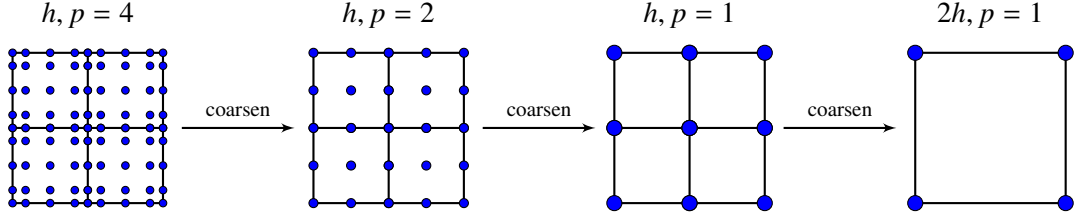
\begin{figure}
    \centering
    \begin{tikzpicture}
	% p=4
	\draw[thick] (0,0) grid +(2,2);
	\foreach \e in {0,1}
	\foreach \f in {0,1}
	\foreach \x in {0,0.1727,0.5,0.8273, 1.0}
	\foreach \y in {0,0.1727,0.5,0.8273, 1.0} 
		\draw[fill=blue] (\e+\x, \f+\y) circle (0.05);

	\node at (1,2.5) {\small $h,p=4$};
	
	\draw[-latex',thick] (2.25, 1) -- node[above] {{\scriptsize coarsen}} (3.75, 1);

    % p=2
    \draw[thick] (4,0) grid +(2,2);
	\foreach \e in {0,1}
	\foreach \f in {0,1}
	\foreach \x in {0,0.5,1.0}
	\foreach \y in {0,0.5,1.0} 
		\draw[fill=blue] (4+\e+\x, \f+\y) circle (0.06);

	\node at (5,2.5) {\small $h,p=2$};
	
	\draw[-latex',thick] (6.25, 1) -- node[above] {{\scriptsize coarsen}} (7.75, 1);
    
    % p=1
    \draw[thick] (8,0) grid +(2,2);
	\foreach \e in {0,1}
	\foreach \f in {0,1}
	\foreach \x in {0,1.0}
	\foreach \y in {0,1.0} 
		\draw[fill=blue] (8+\e+\x, \f+\y) circle (0.1);

	\node at (9,2.5) {\small $h,p=1$};
	
	\draw[-latex',thick] (10.25, 1) -- node[above] {{\scriptsize coarsen}} (11.75, 1);
    
    % 2h,p=1
    \draw[thick] (12,0) rectangle +(2,2);
	\foreach \x in {0,2.0}
	\foreach \y in {0,2.0} 
		\draw[fill=blue] (12+\x, \y) circle (0.1);

	\node at (13,2.5) {\small $2h,p=1$};

	\end{tikzpicture}
    
    \caption{Illustration of the coarsening strategy of GIAMG for quadrilateral elements with nodal basis of polynomial order $p=4$. First, we $p$-coarsen from $p=4$ to $p=2$, and then again $p$-coarsen from $p=2$ to $p=1$. The total number of elements---and therefore the mesh---does not change during $p$-coarsening. Once we reach linear elements, we perform $h$-coarsening to obtain a coarser mesh and smaller system.}
    \label{fig:p_illu}
\end{figure}

There are two major possibilities for constructing a multigrid hierarchy for high-order discretizations: (1) algebraical $h$-coarsening, in which the mesh is coarsened based on the specific matrix properties; and (2) $p$-coarsening, in which the problem is coarsened by reducing the polynomial order of the basis functions and degrees of freedom within each element, resulting in a local coarsening.
While using the geometric information within AMG is not really that beneficial for traditional $h$-coarsening (both problems are equivalent to graph coarsening), $p$-coarsening can be applied in an element-local fashion, and importantly helps sparsify the operator before traditional AMG coarsening is applied. The overall coarsening strategy is illustrated in Figure~\ref{fig:p_illu}.
In this work, we thoroughly evaluate this concept and demonstrate its efficacy in handling complex operators on complex geometries.

For the grids obtained using $p$-coarsening, we effectively increase the sparsity as we go to coarser levels. This enables us to offset the loss of sparsity encountered when generating coarser levels using AMG. This plays a major part in ensuring that the overall scalability remains good.
In addition, if the high-order system is generated using a nodal basis, it interpolates the prolongation matrix values using lower order basis at higher order nodes. If using a modal basis, then the restriction operation is simply an ``injection", i.e., a lower order mode is injected to the coarser level. This indicates the prolongation and restriction operators will only have zero- or unity-valued entries, and the sparsity of coarser level matrices is further reduced compared to using a nodal basis. As a result, it is extremely efficient to perform $p$-coarsening for such high-order systems. %The comparison between nodal and modal basis p coarsening is shown in Figure~\ref{fig:modal_vs_nodal} (\textcolor{red}{need to remake}). 

%Note when using modal basis, 
%In this work, we focus on high-order systems generated from modal basis using {\em Nektar++} -- a well-developed open-source hp element library.      
% \begin{figure}
%     \centering
%     \includegraphics[width=0.7\linewidth]{figures/m_n_comp.jpg}
%     \caption{need to remake}
%     \label{fig:modal_vs_nodal}
% \end{figure}

%review some related work
\par $\boldsymbol{Related\,\, work}$: Multigrid for high-order or spectral/$hp$ elements has been studied as early as in the 1980s. 
%In [4], the authors observe that point smoothers such as the simple Jacobi method result in resolution-independent convergence rates for higher order elements on simple one-dimensional and two-dimensional geometries. Initial theoretical evidence for this behavior is given in [5], where multigrid convergence is studied for one-dimensional spectral methods and spectral element problems. 
The use of $p$-coarsening is rather common in the context of high-order discontinuous Galerkin discretizations~\cite{fidkowski2005p,helenbrook2006application}, but it has also been used for continuous finite element discretizations~\cite{helenbrook2003analysis,huismann2019scaling} and mixed finite element method~\cite{may2015scalable}. A popular strategy for high-order discretizations on unstructured meshes, for which $p$-coarsening is challenging, is to assemble a low-order approximation of the high-order system and use an algebraic multigrid method to invert the low-order (and thus much sparser) operator~\cite{brown2010efficient,kim2007piecewise,olson2007algebraic,deville1990finite,canuto2010finite}. In~\cite{heys2005algebraic}, this approach is compared with the direct application of algebraic multigrid to the
high-order operator, and the authors find that one of the main difficulties is the assembly of the 
%high-order matrices 
low-order approximation required by algebraic multigrid methods. To the best of our knowledge, prior works have not considered using geometric information to assist in setting up the grid hierarchy.  

%what we do
\par $\boldsymbol{Contributions}$: While several approaches have been proposed to solve high-order systems using geometric multigrid, or reducing the operator to design AMG-based preconditioners, there has been limited work on applying AMG directly to the high-order systems.  
In addition, there is a distinct lack of 
scalable, open-source solvers for unstructured high-order discretizations targeting large-scale, real-world  applications. 
In this paper, we proposed a geometrically informed algebraic multigrid (GIAMG) for solving high-order finite element systems in engineering applications, especially for the systems resulting from fluid mechanics as our target application. In a broad sense, the GIAMG method would work for physics problems that can be cast to symmetric systems (operators such as Poison, Helmholtz, etc), depending on the properties of matrices. We also developed a high-performance computing code making use of open-source finite element and AMG libraries. This GIAMG solver employs a $p$-coarsening strategy with assistance of minimal provided geometric information followed by an AMG $h$-coarsening after the order of basis functions that form the original high-order systems is reduced to linear (or some low order). We perform comprehensive experiments 
%by solving a 3D helmholtz operator and a 3D incompressible Navier--Stokes operator, 
and comparisons with other multigrid methods to demonstrate the performance of this GIAMG solver.    

%
%some notes of our work
While this work is motivated by developing a scalable parallel solver for high-order finite element systems, all the matrices used in the experiments are formed 
%using modal basis, specifically 
from a well-recognized open-source spectral/$hp$ element library -- {\em Nektar++} \cite{Cantwell2015,MOXEY2020107110}. %We briefly discuss the difference of modal basis and nodal basis in the paper, but the logic of this GIAMG method persists with minor modifications in the implementation. 
This work is also based on our previously developed smoothed-aggregate AMG (SAAMG) library -- Saena~\cite{saena} for the $h$-coarsening in the GIAMG approach. 
The code for the GIAMG solver is available on \url{https://github.com/paralab/saena_giamg} under the MIT license. 
%The code for the GIAMG solver will be released under the MIT license. 
%
While we have used the open-source {\em Nektar++} and Saena libraries, the GIAMG code can easily be modified to use other high-order finite-element packages such as {\texttt deal.ii} \cite{dealii},  MOOSE~\cite{lindsay2022moose} and MFEM~\cite{mfem}, and AMG packages such as ML \cite{amg-ml} and Hypre \cite{hypre}.

\par $\boldsymbol{Limitations}$: Limited geometric information, which is available in {\em Nektar++} as well as most high-order finite element packages, currently needs to be passed to the GIAMG solver. Note, the GIAMG solver is not limited to be used with {\em Nektar++}, we only use {\em Nektar++} as an example demonstration in this work. The GIAMG solver has an interface to take necessary geometric information passed from high-order finite element packages and convert the information to its own data structure. The information needed is limited to identifying which degrees-of-freedom (dof) can be grouped together as an element, as well as their role within the elemental basis. This makes GIAMG solver not a pure AMG solver (as a result, so-called geometrically informed). It is challenging to automatically identify such information purely from the input matrices, which is the bottle-neck to make GIAMG a complete \textit{black-box} AMG solver for solving high-order finite element systems. We are currently working on efficient methods to identify this information directly from the discretized operators. In addition, the GIAMG solver currently only works with symmetric systems for optimal convergence, and therefore treatment is required for non-symmetric terms in Navier--Stokes equations (e.g., using the velocity correction scheme).   

The rest of the paper is organized as follows. 
%In Section~\ref{theory}, we describe the theory of $p$-coarsening and the GIAMG method. 
In \S~\ref{implem}, we describe in detail the algorithms and implementation of GIAMG method used for high-order systems. In \S~\ref{results}, we present comprehensive experiments and comparisons with other available multigrid packages using 3D test problems. Finally, in \S~\ref{conc} we draw conclusions and discuss future research directions.

% \section{Theory and Method 
% }\label{theory}
% \subsection{Benefit of $p$-coarsening}\label{p_benefit}

\section{Algorithms and implementations}\label{implem}

Consider a linear system 
\begin{align}
Ax=b
\end{align}
where $A$ is a linear operator, $b$ is a right-hand-side vector, and $x$ is the solution vector. The multigrid v-cycle is the standard building block for both multigrid as a solver, as well as when used as a preconditioner with the conjugate gradient method. We therefore start by summarizing the  overall process of applying a v-cycle of the proposed GIAMG method. The v-cycle algorithm is presented in Algorithm~\ref{GIAMG}.

\begin{algorithm}[ht!]
  \caption{GIAMG process}\label{GIAMG}
  \begin{algorithmic}[1]
  %\State {loop each element}
    \Procedure{v-cycle}{}
        \If {coarsest level}
            \State Directly solve the coarsest matrix
            \State return
        \EndIf
        \State Presmooth
        \State Compute residual
        \If{$p$-coarsening level}
            \State $p$-restrict residual
        \Else 
            \State $h$-restrict residual
        \EndIf
        \State {\footnotesize{V-CYCLE}}
        \If{$p$-coarsening level}
            \State $p$-interpolate error
        \Else 
            \State $h$-interpolate error
        \EndIf
        \State Correct solution
        \State Postsmooth
    \EndProcedure
  \end{algorithmic}
\end{algorithm}
At each level, we employ a Jacobi-accelerated Chebyshev smoother for presmoothing. Prolongation and restriction operators are formed according to whether the current grid employs $p$-coarsening or SAAMG coarsening. The coarser level operator is computed using a Galerkin projection as $A_c = RA_fP$, where $A_f$ is the finer level operator, $P$ and $R$ are prolongation and restriction, respectively, and $P = R^T$. Note that the coarse-grid operator, as well as the prolongation and restriction operators at each level are precomputed and stored during the setup phase. The finer level residual vector is computed after smoothing as $r_f = b_f - A_fx_f$, where $b_f$ is the finer level right-hand-side and $x_f$ is the current estimate of the finer level solution vector, and restricted to the next coarser level as $r_c = Rr_f$, where $r_c$ is coarser level right-hand-side. A coarser level linear system (error equation) is formed as $A_c e_c=r_c$, where $e_c$ is the coarse-level error vector. This process is recursively repeated until the coarsest level, where we employ a direct solver (superLU \cite{superLU}) to solve the coarsest system. 
The coarse-level error vector is then prolonged to be the finer level error vector as $e_f=Pe_c$, and the finer level solution vector is corrected as $x_f=x_f+e_f$. A postsmoothing (similar as presmoothing) is applied to the corrected solution vector at each level. This process is again recursively repeated until it reaches the finest level to obtain the final solution.
The GIAMG approach is used as a preconditioner in this work since the multigrid method is more robust when used as a preconditioner in conjunction with CG~\cite{sundarNLA15}. Since we are dealing with symmetric high-order systems, a GIAMG preconditioned conjugate gradient (PCG) iterative solver is employed. One v-cycle of GIAMG is called during the preconditioning, and the rest is the same as standard PCG method.  
Note during the $p$-coarsening, we have the flexibility to choose how aggressive, i.e., how many orders of basis to reduce at each level, which can affect the overall convergence (with the SAAMG coarsening levels fixed). We will demonstrate the convergence performance of different $p$-coarsening strategies in \S~\ref{NS}. In addition, the computationally dominant part of the GIAMG method is the pre- and postsmoothing during the  v-cycle. 
%(probably at the first level). 
%It should in theory cost much more than the summation of all subsequent coarser levels, and dominate the overall cost. 
Thus, the performance of the Chebyshev smoother is important, and we want as few iterations (per v-cycle) as possible while still achieving good smoothing at each level. We also demonstrate the effect of smoothing on the overall convergence and total solve time in \S~\ref{NS}. 

One challenge with the GIAMG apporach is that additional geometric information, which should be provided by the finite element assembly library, is required at the finest level and needs to be reconstructed at coarser levels to inform how to pick dofs that remain at coarser levels during $p$-coarsening. This information can be discarded after all the coarse grids are constructed.
%(i.e., coarser levels with reduced order of basis functions). 
Note, ideally, it would be optimal to automatically identify such information purely from the matrices. However, it is not easy to extract this information in an efficient and scalable fashion.

\subsection{Required geometric information}
%Therefore, in this work, we propose a geometrically-informed algebraic multigrid (GIAMG) for solving high-order finite element systems.
%
Since we reduce the basis order within each element, we need to know the elemental dof collection for each element, and also the ordering of the elemental dofs, which implies the basis order at each dof so that we know how to select the dofs that are retained at the next coarser level during $p$-coarsening. In this work, we use {\em Nektar++} to assemble the high-order system, and an example of how {\em Nektar++} orders elemental dofs based on its basis definition is described in the Appendix. In GIAMG solver, it defines an \textit{l2g} map, which is the local to global mapping. 
%The input of this l2g map is the elemental dof index, and the output is the index in the global system for this dof. 
It is a 2D vector in which a row represents each element, and a column represents collection of elemental dofs. The value of each entry in the \textit{l2g} map is the dof index in the process matrix. The ordering of entry values in each column (i.e., the ordering of the elemental dofs) in this map implies the basis orders. 
This procedure is completely parallel, and can be performed independently on each process without any communication. In GIAMG solver, it has an application program interface (API) to take the local to global map from high-order finite element libraries and convert the map to its own \textit{l2g} map data structure, which requires minimal modification specific to the finite element library data structure. In the parallel case~\footnote{In the parallel case, we use ``local to global" to represent elemental to process mapping, and ``global to universal" to represent process to universal mapping.}, GIAMG solver defines a \textit{g2u} map from process dofs to universal system indices. The \textit{g2u} map is a 1D array with the indices representing dof indices in each process (values of \textit{l2g} map) and values representing the dof indices in the universal system. GIAMG solver also has an API to convert the data structure of this global to universal map. The local to global and global to universal maps, which are generally available in most of high-order finite element libraries, are the only external arguments that are required in the GIAMG solver.
%Using these two maps, we can identify which dofs are "dropping" to the next lower dimensional space. 

%A map from the local elemental dofs to the global matrix indices is defined in {\em Nektar++} and passed to the GIAMG solver.  
%

%~\footnote{We are currently able to identify elemental dof collections, or cliques, purely from an input matrix. However, to identify the elemental dof ordering still remains challenging to us.}. 
%For example, using {\em Nektar++}, we need to provide the following geometry information to the GIAMG solver. 

\subsection{Computation of \textit{l2g} map at each level}
%{\em Nektar++} defines an \textit{l2g} map, which is the local to global mapping. 
%GIAMG solver defines an \textit{l2g} map, which is the local to global mapping. 
%The input of this l2g map is the elemental dof index, and the output is the index in the global system for this dof. 
%It is a 2D vector in which a row represents each element, and a column represents collection of elemental dofs with the ordering defined by Equation~\ref{modal_basis} in each element. The value of each entry in the \textit{l2g} map is the dof index in the process.
%The ordering of the columns follows the pattern we discussed previously in Section~\ref{elemental_ordering}. 
For a {\em Nektar++} hexahedron element type, to extract the dofs that remain to the next coarser level, for example, if we perform a $p \rightarrow p/2$ coarsening (assuming $p$ is an even number), and denote the finer level and coarser level elemental dof collections in each element as $df_{arr}$ and $dc_{arr}$, respectively, we have Algorithm~\ref{extract}. 
\begin{algorithm}[ht!]
  \caption{Extract elemental dofs for next level}\label{extract}
  \begin{algorithmic}[1]
  %\State {loop each element}
    \Procedure{extract\_next\_dof}{$df_{arr}$}
    %\State create $R$ and {\color{red}$n_R$\Comment{$n_R$ is the bin number}}
        \For{$i=0:p/2$}
            \For{$j=0:p/2$}
                \For{$k=0:p/2$}
                    \State Append $df_{arr}((p+1)^2\times i+(p+1)\times j+k)$ to $dc_{arr}$
                \EndFor
            \EndFor
        \EndFor
        \Return $dc_{arr}$
    \EndProcedure
  \end{algorithmic}
\end{algorithm}
For other element types, we similarly follow the elemental dof ordering implied by the definition of basis in those elements (similar as in the Appendix). However, triangular-shaped (including prism) elements have more complicated orderings. 
%Note, the dof ordering in the $dc_{arr}$ also matches the ordering in the definition of modal basis at the coarser level.  

With the \textit{l2g} map and extraction of dofs in the next level, we can first collect all dofs that remain at the next level in an array, denoted as \textit{collect}, by looping over the rows of the \textit{l2g} map (i.e., looping over elements), as illustrated in Algorithm~\ref{get_collect}.  
\begin{algorithm}[ht!]
  \caption{Collect all dofs in next level}\label{get_collect}
  \begin{algorithmic}[1]
  %\State {loop each element}
    \Procedure{dof\_next}{$l2g$}\Comment{}
    %\State create $R$ and {\color{red}$n_R$\Comment{$n_R$ is the bin number}}
        \For{$e_i = 1:\text{size}(l2g,1)$}
            \State $dc_{arr} = \text{\scriptsize{EXTRACT\_NEXT\_DOF}}(l2g(e_i, :))$
            \For{$dof = 1: \text{length}(dc_{arr})$}
                \If {$dc_{arr}(dof) \in collect$}
                    \State Skip
                \Else
                    \State Append to collect
                \EndIf
            \EndFor
        \EndFor
        \State Sort $collect$ 
        \Return $collect$
    \EndProcedure
  \end{algorithmic}
\end{algorithm}
The indices of the $collect$ will be the values of the next level local to global map, denoted as $l2g_{next}$. Then we build a Hash map, denoted as $hash_{l2g}$, mapping from the $collect$ values (representing finer level $l2g$ values) to its indices (representing coarser level $l2g_{next}$ values). Note, since we perform a sort for $collect$, the coarser level system will have a totally different permutation, which is computed during the sort, compared with the finer level. However, this does not matter since even though we change the permutation in the coarser levels during residual restriction, when we prolong the error back to finer levels, it will be projected back to the original permutation at each level. The $l2g_{next}$ can be then computed as in Algorithm~\ref{compute_l2g_next}. The computation of $l2g_{next}$ will need to be performed at each level. Note, since $p$-coarsening only reduce $p$ order within elements, the total number of elements remains the same in each level, i.e., the total number of rows is the same in $l2g$ and $l2g_{next}$. Also, in the computation of $l2g_{next}$, boundary nodes need to be carefully handled. 
\begin{algorithm}[ht!]
  \caption{Compute $l2g_{next}$ map in the next level}\label{compute_l2g_next}
  \begin{algorithmic}[1]
  %\State {loop each element}
    \Procedure{l2g\_next}{$l2g$}
        \For{$e_i = 1:\text{size}(l2g,1)$}
            \State $dc_{arr} = \text{\scriptsize{EXTRACT\_NEXT\_DOF}}(l2g(e_i, :))$
            \For{$dof = 1: \text{length}(dc_{arr})$}
                \State $l2g_{next}(e_i, dof) = hash_{l2g}(dc_{arr}(dof))$
            \EndFor
        \EndFor
    \EndProcedure
  \end{algorithmic}
\end{algorithm}

\subsection{Parallel implementation}
In addition to the $l2g$ map, we need another $g2u$ map for parallel computing. In essence, we have three numberings, the local or elemental numbering, the global or process-local numbering and the universal or the numbering across all dofs on all processes.
%The input of this g2u map is the value of l2g map for each dof, and the output is the index in the universal system for this dof. 
The $g2u$ map is also read in at the finest level, and we need to re-compute it at each coarser level. We first collect coarser level dofs in each process, and synchronize between all processes to obtain $collect_{univ}$, as shown in Algorithm~\ref{sync_g2u}.
\begin{algorithm}[ht!]
  \caption{Obtain universal $collect_{univ}$}\label{sync_g2u}
  \begin{algorithmic}[1]
  %\State {loop each element}
    \Procedure{dof\_next\_univ}{$l2g$}
        \State $collect$ = {\scriptsize{DOF\_NEXT}}($l2g$)
        \State Gather $collect$ to $collect_{univ}$
        \State Sort and remove duplicates in $collect_{univ}$
    \EndProcedure
  \end{algorithmic}
\end{algorithm}
We then construct a Hash map, $hash_{g2u}$, mapping from $collect_{univ}$'s values to its indices, similarly as $hash_{l2g}$. Now we are able to construct the next coarser level $g2u$ map, $g2u_{next}$, by the auxiliary of $collect$, $g2u$ and $l2g_{next}$, as shown in Algorithm~\ref{compute_g2u_next}. We loop over rows of $l2g_{next}$ (i.e., elements in each process), and obtain each dof's $l2g_{next}$ value, which corresponds to the index of the $collect$ array. 
%mapping to finer level $l2g$ value for this dof. 
This $collect$ value further corresponds to the index of the finer level $g2u$ array.
%mapping to the index in the universal system for this dof. 
Using the $hash_{g2u}$, we finally compute the value in the $g2u_{next}$ array for this dof.  
\begin{algorithm}[ht!]
  \caption{Compute $g2u_{next}$ map in next level}\label{compute_g2u_next}
  \begin{algorithmic}[1]
    \State $collect$ = {\scriptsize{DOF\_NEXT}}($l2g$)    \Procedure{g2u\_next}{$l2g_{next}, g2u$}
        \For{$e_i = 1:\text{size}(l2g_{next},1)$}
            \For{$dof = 1:\text{length}(l2g_{next}(e_i,:))$}
                \State $coarse_{glb} = l2g_{next}(e_i,dof)$ 
                \State $fine_{glb} = collect(coarse_{glb})$;
                \State $fine_{univ} = g2u(fine_{glb})$;
                \State $g2u_{next}(coarse_{glb}) = hash_{g2u}(fine_{univ})$;
            \EndFor
        \EndFor
    \EndProcedure
  \end{algorithmic}
\end{algorithm}

\subsection{Forming prolongation}

With the help of $l2g$ and $g2u$ maps at each level, we are able to compute prolongation $P$ as shown in Algorithm~\ref{compute_P}. We construct a sparse matrix  in the {\footnotesize{COO}} format \footnote{List of tuples of row, column, value.} for $P$, denoted as $P_{coo}$.  We loop over elements in each process, and loop over each coarser level dof within this element. We find values in both coarser level $g2u_{next}$ map and finer level $g2u$ map for this dof, and set them to be the row and column indices in the universal system, respectively, and assign 1.0 (``injection" for modal basis~\footnote{GIAMG solver also has implementation for nodal basis, and it constructs prolongation matrix using interpolation instead of ``injection".} ) as the value in the {\footnotesize{COO}} structure. We append each {\footnotesize{COO}} structure of coarser level dofs in $P_{COO}$ and remove duplicates.  
\begin{algorithm}[ht!]
  \caption{Compute prolongation $P$}\label{compute_P}
  \begin{algorithmic}[1]
  %\State {loop each element}
    \Procedure{comp\_prol}{$l2g, l2g_{next}, g2u, g2u_{next}$}
        \For{$e_i = 1:\text{size}(l2g,1)$}
            \State $dc_{arr}$ = {\scriptsize{EXTRACT\_NEXT\_DOF}}($l2g(e_i, :)$)
            \For{$dof = 1: \text{length}(dc_{arr})$}
                \State $P_{col}$ = $g2u_{next}(l2g_{next}(e_i,dof))$
                \State $P_{row} = g2u(dc_{arr}(dof))$ 
                \State Append $\left(P_{row}, P_{col}, 1.0\right)$ to $P_{coo}$ 
            \EndFor
        \EndFor
        \State Sort $P_{coo}$ and remove duplicates
    \EndProcedure
  \end{algorithmic}
\end{algorithm}

\subsection{h-coarsening}
After reducing the basis order to the linear level, we employ the standard smoothed-aggregation AMG (SAAMG) approach to create coarser grids for h-coarsening. In this work, our previously developed open-source SAAMG library, Saena \cite{saena}, is employed for this purpose. The aggregates and corresponding tentative prolongation are generated using a minimal independent set (MIS) approach. To smooth the tentative prolongation, an SPAI \cite{broker2002sparse} smoother is employed.  
For more details about Saena, we refer readers to~\cite{saena}.

\section{Experiments and results}~\label{results}
In this section, we evaluate the performance of the GIAMG solver using the Helmholtz and the incompressible flow problems. 

\subsection{Helmholtz operator}\label{helmholtz}
We present a 3D Helmholtz operator posed on a unit cube domain with 8 elements using spectral/$hp$ element. The equation is written as
\begin{align}
     - \nabla^2 \text{u} + \lambda \text{u} = \text{f}.
\end{align}
%
% which can be cast into a weak form as
% \begin{align}
%      \int \nabla w \cdot \nabla u + \int w \lambda u = \int w f
% \end{align}
In this case, we manufacture the forcing term $\text{f}$ using an analytical solution as below
\begin{align}
    \text{f} &= -(\lambda +3\pi^2)\text{sin}(\pi x)\text{sin}(\pi y)\text{sin}(\pi z)\\
    \text{u} &= \text{sin}(\pi x)\text{sin}(\pi y)\text{sin}(\pi z)
\end{align}
$\lambda$ is chosen to be 1.0~\footnote{We note that the results do not extend to high-frequency Helmholtz problems.}. Zero Dirichlet boundary condition is specified at all boundaries, and zero initial guess is employed. We only show the convergence results in this section, and we present detailed profiling and scaling results of the GIAMG approach in \S~\ref{NS}.

\subsubsection{Convergence results}
\begin{figure}
\centering
\begin{tikzpicture}
    \begin{semilogyaxis}[
    xlabel={Number of iterations},ylabel={Relative residual},xmin=0,ymin=1e-14,ymax=10,
    width=3in,
    height=2.5in,
    ymajorgrids = true,
    legend entries={p=3,p=4,p=5,p=6,p=7,p=8,p=9,p=10},legend style={nodes={scale=0.8, transform shape}},legend pos=north east, legend columns=3]
        \addplot [very thick,smooth,color=blue] table[] {helmholdz/p3_paper.txt};
        \addplot [very thick,smooth,color=green] table[] {helmholdz/p4_paper.txt};
        \addplot [very thick,smooth,color=red] table[] {helmholdz/p5_paper.txt};
        \addplot [very thick,smooth,color=black] table[] {helmholdz/p6_paper.txt};
        \addplot [very thick,smooth,color=yellow] table[] {helmholdz/p7_paper.txt};
        \addplot [very thick,smooth,color=brown] table[] {helmholdz/p8_paper.txt};
        \addplot [very thick,smooth,color=magenta] table[] {helmholdz/p9_paper.txt};
        \addplot [very thick,smooth,color=cyan] table[] {helmholdz/p10_paper.txt};
    \end{semilogyaxis}
\end{tikzpicture}
\caption{Convergence results of GIAMG from $p$ = 3 to 10 for the Helmholtz operator. 
%\textcolor{blue}{[SX] should we only plot before 1e-6 to show it is order-independent?}
} \label{fig:saena_p}
\end{figure}
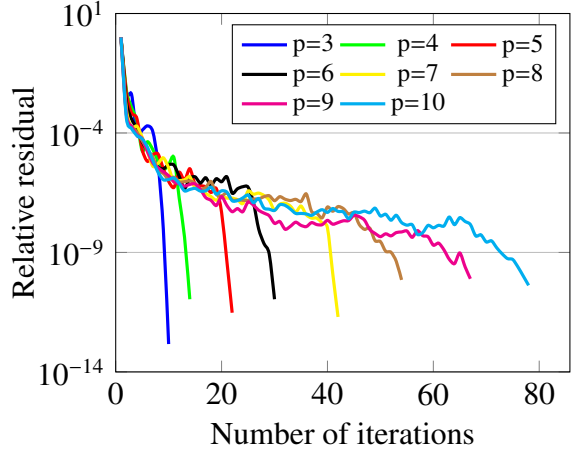

We first plot the convergence for basis order from $p = 3$ to $10$ using the GIAMG method by reducing the basis order by 2 at each level as shown in Figure.~\ref{fig:saena_p}. Note, in these cases, we only need to perform $p$-coarsening.
% \begin{figure}
%     \centering
%     \includegraphics[width=0.7\linewidth]{figures/saena_p_more.eps}
%     \caption{Convergence results of present work from p = 3 to 10 \textcolor{blue}{[SX] should we only plot before 1e-6 to show it is order-independent?}}
%     \label{fig:saena_p}
% \end{figure}
We set the convergence criterion as the relative residual being less than $10^{-10}$. %\textcolor{blue}{[SX] or maybe 1e-6 since the convergence is more order-independent before 1e-6}. 
We can see as we increase the basis order, the number of iterations does not significantly increase before a relative residual around $10^{-6}$. The increase in the overall number of iterations is primarily because of the inability of the Chebyshev smoother to handle high-order discretizations, as has been shown in prior work \cite{sundarNLA15}.

% implying the GIAMG approach is order-independent till a reasonably small tolerance for convergence.

\subsubsection{Comparison with other solvers at p = 8}\label{helmholtz_p8}
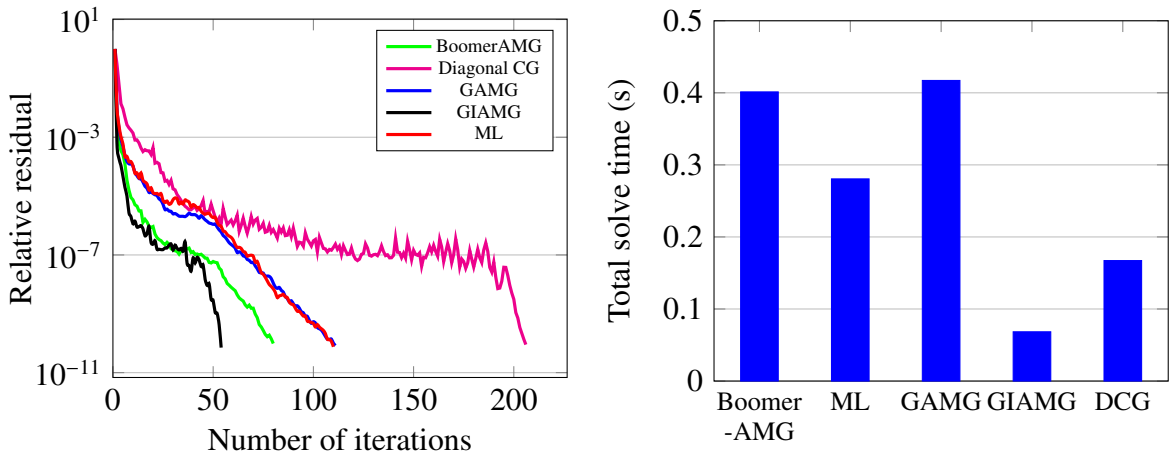
\begin{figure}[ht]
  \quad
  \begin{subfigure}[b]{0.45\linewidth}
    \begin{tikzpicture} %[scale=0.85]
      \begin{semilogyaxis}[width=3in,height=2.5in,
      xlabel={Number of iterations},xmin=0,ylabel={Relative residual},
      ymajorgrids = true,
      legend entries={BoomerAMG,Diagonal CG,GAMG,GIAMG,ML},legend style={nodes={scale=0.6, transform shape}},legend pos=north east]
        \addplot [very thick, color=green] table[] {helmholdz/boomeramg_1e-10_paper.txt};
        \addplot [very thick, color=magenta] table[] {helmholdz/dcg_1e-10_paper.txt};
        \addplot [very thick, smooth, color=blue] table[] {helmholdz/gamg_1e-10_paper.txt};
        \addplot [very thick, smooth,color=black] table[] {helmholdz/giamg_1e-10_paper.txt};
        \addplot [very thick,smooth,color=red] table[] {helmholdz/ml_1e-10_paper.txt};
      \end{semilogyaxis}
    \end{tikzpicture}
    %\caption{Convergence} 
    \label{fig:P8_comp_converge}
  \end{subfigure}
  \quad
  \begin{subfigure}[b]{0.45\linewidth}%\centering
    \begin{tikzpicture}[baseline = {(0,-1.06)}]
    \begin{axis}[
        width=3in,height=2.51in,
        %height=2.8in,
        %scale only axis,
        %clip=false,
        %separate axis lines,
        %axis on top,
        ymajorgrids = true,
        xmin=0.5,
        xmax=5.5,
        xtick={1,2,3,4,5},
        %x tick style={draw=none},
        x tick label style={%rotate=45, 
        align=center},
        xticklabels={\footnotesize Boomer \\[-3pt] \footnotesize {-AMG},\footnotesize ML,\footnotesize GAMG,\footnotesize GIAMG,\footnotesize DCG},
        ytick={0,0.1,0.2,0.3,0.4,0.5},
        ymin=0,
        ymax=0.5,
        ylabel={Total solve time (s)},
        every axis plot/.append style={
          ybar=2pt,
          %x=1cm,
          bar width=0.5cm,
          bar shift=0pt,
          fill
        }
      ]
      \addplot[blue]coordinates {(1,4.1527e-01 - 1.3987e-02)};
      \addplot[blue]coordinates{(2,3.1273e-01 - 3.2169e-02)};
      \addplot[blue]coordinates{(3,5.7788e-01 - 1.6072e-01)};
      \addplot[blue]coordinates{(4,0.00129*53)};
      \addplot[blue]coordinates{(5,8.196394e-04*204)};
    \end{axis}
  \end{tikzpicture}
    %\caption{Solve time} 
    \label{fig:P8_comp_timing}
  \end{subfigure}
  \caption{Convergence and solve time comparisons of GIAMG with PETSc GAMG, ML, BoomerAMG and DCG at $p$ = 8 for the Helmholtz operator.}
  \label{fig:P8_comp}
\end{figure}

We  compare the performance of the GIAMG method with other well-recognized multigrid packages, including PETSc GAMG, ML, and Hypre BoomerAMG solvers, at $p = 8$. Note, we call ML and BoomerAMG also via their PETSc interface. Since there are many components that affect the performance of multigrid solvers, and varied approaches have their own strategies and parameters to handle different parts of multigrid, it is not possible to perform an exact comparison among different multigrid approaches. As a result, to initiate a reasonable comparison between all approaches, we keep some key parts in multigrid relatively comparable by tuning parameters in these parts for each approach according to our best knowledge. First of all, we use v-cycle for all cases. Secondly, note that in the coarsening part, GIAMG uses $p$-coarsening, while PETSc GAMG and ML use smoothed-aggregate method, and BoomerAMG uses the classical C-F splitting. The way we try to make the multigrid hierarchy fairly comparable is to ensure they all have the same number of levels and at the coarsest level, they have relatively comparable matrix sizes and numbers of nonzeros.
%by tuning coarsening schemes and parameters in these approaches. 
Finally, we ensure that they all use the Chebyshev smoother during pre- and postsmoothing with the same number of smoothing iterations. The breakdown of the key setup of different approaches 
%that affect the performance of multigrid solver 
is shown in Table~\ref{tab:setup}.
\begin{table}
    \centering
    \begin{tabular}{|p{23mm}|p{25mm}|p{24mm}|p{18mm}|p{25mm}|}
    \hline
        & Smoother /\newline Iterations & Coarsening \newline scheme & Hierarchy \newline Levels & Coarsest level \newline size/nonzeros\\
    \hline
    GAMG & Chebyshev / 3& Default SA & 4 & 21/441\\
    \hline
    ML & Chebyshev / 3 & Uncoupled&4 & 23/365\\
    \hline
    BoomerAMG &Chebyshev / 3 &HMIS & 4 &26/310\\
    \hline
    GIAMG&Chebyshev / 3 &$p$-coarsening &4 &27/343 \\
    \hline
    \end{tabular}
    \caption{Key setup used in each multigrid approach at $p = 8$ for the Helmholtz operator.}
    \label{tab:setup}
\end{table}
%
%The detailed Petsc options used in each approach is shown in Appendix. 
We also provide the performance of a diagonally preconditioned CG solver (DCG) implemented in {\em Nektar++} for comparison. Since the problem size is not too large, we can use a single core to test all the methods. The comparisons of the convergence and total solve time at $p = 8$ are plotted in Figure~\ref{fig:P8_comp}.
% \begin{figure}
%     \centering
%     \includegraphics[width=0.45\linewidth]{figures/P8_comp_more.eps}
%     \includegraphics[width=0.45\linewidth]{figures/P8_comp_timing.eps}
%     \caption{Convergence comparison of present work with PETSc GAMG, ML, BoomerAMG and DCG at p = 8}
%     \label{fig:P8_comp}
% \end{figure}
%
% \begin{figure}
%     \centering
%     \includegraphics[width=0.7\linewidth]{figures/P8_comp_timing.jpg}
%     \caption{Comparison of solve time (without considering multigrid setup time) between present work, PETSc GAMG, ML, BoomerAMG and DCG at p = 8}
%     \label{fig:p8_comp_timing}
% \end{figure}
We can see our GIAMG's convergence and total solve time are both the best among all methods. 
% Interestingly, while all AMG approaches besides GIAMG took longer to solver compared with DCG, in spite of DCG taking more iterations overall. 
Interestingly, in spite of DCG taking more iterations overall for convergence, it costs less solve time compared with other AMG approaches, indicating AMG approaches (including GIAMG) take longer time to solve for one PCG iteration compared with DCG.
GIAMG is the only AMG solver that is faster than DCG in terms of total solve time. Note that the overall solve time for GIAMG is significantly better than the other AMG solvers in part due to the fewer iterations it requires, but also because of the lower cost per iteration due to increased sparsity of the coarse-grid operators.

\subsection{Incompressible flow problem}\label{NS}
For the next evaluation, we solve for the 3D incompressible flow problem with a geometry based upon flow simulations with application to thermal transport in nuclear fusion modeling. The incompressible Navier--Stokes equations may be written as 
\begin{align}
    \frac{\partial \textbf{u}}{\partial t}+ \textbf{u} \cdot \nabla \textbf{u}  &= - \nabla \text{p} + \nu \nabla^2 \textbf{u} +  \textbf{f}
\\
\nabla \cdot \textbf{u}  &= 0
\end{align}
where $\textbf{u}$ is the velocity, $\text{p}$ is the pressure and $\nu$ is the kinematic viscosity. The governing equations are solved using a widely-adopted stiffly stable velocity correction scheme which splits the pressure and velocity fields. The velocity and pressure systems are decoupled resulting in four systems of rank $N$ instead of a single coupled system of rank $4N$, %in the incompressible Navier--Stokes equations, 
where $N$ denotes the total number of discretization points. The velocity correction scheme integrates the non-linear convection (i.e. first derivative) term explicitly in time whilst the pressure equation and diffusion (i.e. second derivative) term are integrated implicitly, requiring the solution of the Poisson and Helmholtz operators, respectively. Thus, the decoupled matrices are elliptic and symmetric.
%Splitting schemes are typically favourite for their numerical efficiency since for a Newtonian fluid the velocity and pressure are handled independently, requiring the solution of three (in two dimensions) elliptic systems of rank N (opposed to a single system of rank 3N solved in the Stokes problem). 
%However, a drawback of this approach is the splitting scheme error which is introduced when decoupling the pressure and the velocity system, although this can be made consistent with the overall temporal accuracy of the scheme by appropriate discretisation of the pressure boundary conditions. 
%When the incompressible Navier-Stokes equations are solved using the velocity correction splitting scheme, a stiffly-stable time integration is applied. 
%(derived from the work of Karniadakis, Israeli and Orszag). 
For detailed information on the velocity correction scheme implemented in {\em Nektar++}, we refer readers to~\cite{KARNIADAKIS1991414,MOXEY2020107110}. We note that solving the pressure field is typically the most challenging part in the velocity correction scheme, and therefore we focus on solving for pressure in the numerical experiments below.

\begin{figure}
    \centering
    \includegraphics[width=0.31\linewidth]{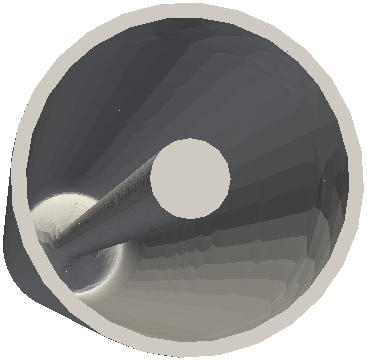}
    \hspace{2cm}
    \includegraphics[width=0.3\linewidth]{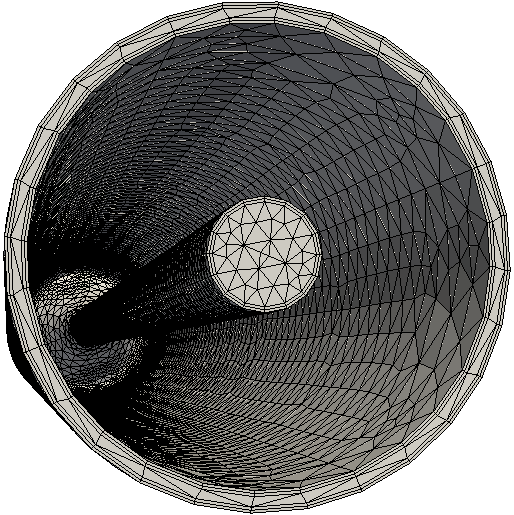}
    \vspace{5mm}
    \caption{Geometry and mesh of the computational domain.}
    \label{fig:nano_geo_mesh}
\end{figure}

The computational domain and mesh are shown in Figure~\ref{fig:nano_geo_mesh}, which represents a sample pin-type breeder blanket for a fusion reactor that is used in part to cool the reactor. In this application, a cooling fluid enters through the circular cross-section in the center of the pin, transports heat radiated at the opposite end of the pin and exits through the outer annulus. In this experiment, we do not consider the thermal transport terms, and focus on a under-resolved DNS at a Reynolds number of $Re=10^4$ with a coarse mesh of $90,842$ tetrahedra and $71,696$ prisms, where the latter elements are used to resolve boundary layers in near-wall regions. The coarse nature of these simulations, together with large changes in velocity across a small region at the end of the pin, together with high aspect ratio elements at the walls, makes the systems particularly challenging to precondition at high order.

Boundary conditions are specified as a Poiseuille-type inflow condition supplied on the center of the pin with a standard Neumann outflow condition supplied on the outflow (the outer annulus), and no-slip boundary conditions imposed on the remaining wall surfaces. Initial conditions to assemble the systems are the solutions after solving the governing equations for some time steps and reaching a physically-meaningful state. We only solve for one time step in our experiments since the matrices are not changing with time.
%We use solutions at this physically-meaningful state to assemble the systems. 
Considering that very high order of basis function is typically not favored in practice for such a complex engineering application, we choose a reasonably high order of $p = 5$ as our maximum basis order in this example. The finest level matrix size and number of nonzeros are $6,604,956\times6,604,956$ and $1,129,952,608$, respectively. All the experiments (except the scaling experiment) are performed using 3 nodes with 40 cores on each node, resulting in a total 120 MPI tasks using University of Utah CHPC computing resources.

\subsubsection{Convergence results}\label{NS_converge}
\begin{figure}
\centering
\begin{tikzpicture}
    \begin{semilogyaxis}[width=3.5in,height=2.5in,ymajorgrids = true,ymin=1e-7,ymax=10,
    xlabel={Number of iterations},ylabel={Relative residual},xmin=0,legend entries={p=1,p=2,p=3,p=4,p=5},legend style={nodes={scale=0.8, transform shape}},legend pos=north east]
        \addplot [very thick, smooth,color=red] table[] {nano/P1_new_paper.txt};
        \addplot [very thick, smooth,color=blue] table[] {nano/P2_new_paper.txt};
        \addplot [very thick, smooth,color=magenta] table[] {nano/P3_new_paper.txt};
        \addplot [very thick, smooth,color=green] table[] {nano/P4_new_paper.txt};
        \addplot [very thick, smooth,color=black] table[] {nano/P5_new_paper.txt};
    \end{semilogyaxis}
\end{tikzpicture}
\caption{Convergence results for $p$ from 1 to 5 for the incompressible flow problem.} \label{fig:nano_all_p}
\end{figure}
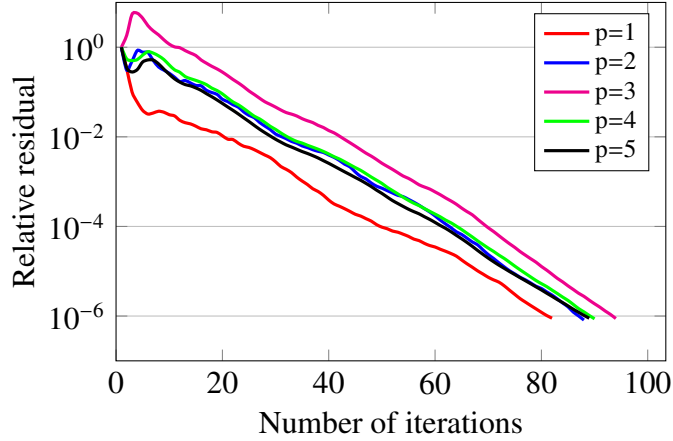

We plot the convergence results using GIAMG for $p$ from 1 to 5 as shown in Figure~\ref{fig:nano_all_p}. For this more complex system, we set the convergence criterion to be the relative residual less than $10^{-6}$. In this case, we only reduce the basis order by 1 at each level during $p$-coarsening followed by 6 additional levels of coarsening using smoothed aggregation. We perform 2 Chebyshev iterations in pre- and postsmoothing. The sizes and numbers of nonzeros of the coarsest level matrices in these cases are slightly different, varying from $303\times303$ to $315\times315$, and $32338$ to $34191$, respectively, which are all reasonably small at the coarsest level for the superLU solver. We can see the number of iterations does not significantly increase as we increase $p$ order. The convergence of $p = 3$ case looks slightly inferior to other cases. However, it still behaves within a reasonable range, and the convergence pattern is similar to that in other cases.
% \begin{figure}
%     \centering
%     \includegraphics[width=0.7\linewidth]{figures/comparison_p_order.eps}
%     \caption{Convergence results for p from 1 to 5}
%     \label{fig:nano_all_p}
% \end{figure}

\subsubsection{Detailed analysis at p = 5}
We present a detailed analysis at $p = 5$ to comprehensively evaluate the performance of our GIAMG solver. All the timing presented in this section is the average across all processes. 6 levels of SA coarsening are used (as in the previous section) in all experiments except for the comparison with other multigrid packages, which is  detailed in \S\ref{petsc_comp}. The same machine is employed (as in the previous section to allow inter-order comparisons) except for the strong scaling experiments, which are presented in \S\ref{scaling}.

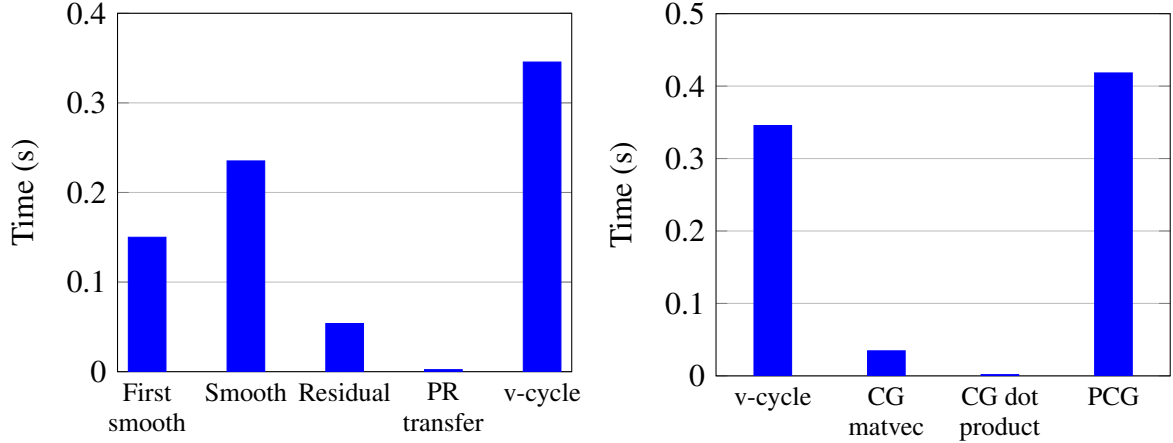
\begin{figure}[ht]
  \quad
  \begin{subfigure}[b]{0.45\linewidth}%\centering
    \begin{tikzpicture}%[scale=0.7]
    \begin{axis}[
        % yticklabel style={
        % /pgf/number format/fixed,
        % /pgf/number format/precision=5
        % },
        %scaled y ticks=false,
        xticklabel style={align=center},
        width=3in,height=2.5in,
        %height=2.8in,
        %scale only axis,
        %clip=false,
        %separate axis lines,
        %axis on top,
        %title={\large{Breakdown time at p = 5 for v-cycle}},
        %title style={yshift=0.75ex},
        ymajorgrids = true,
        xmin=0.7,
        xmax=5.3,
        xtick={1,2,3,4,5},
        x tick style={draw=none},
        xticklabels={\footnotesize First \\[-3pt] \footnotesize smooth,\footnotesize Smooth,
        %\footnotesize superLU,
        \footnotesize Residual,\footnotesize PR \\[-3pt] \footnotesize transfer, \footnotesize v-cycle},
        ymin=0,
        ymax=0.4,
        ytick={0,0.1,0.2,0.3,0.4},
        ylabel={Time (s)},
        every axis plot/.append style={
          ybar=2pt,
          %x=1cm,
          bar width=0.5cm,
          bar shift=0pt,
          fill
        }
      ]
      \addplot[blue]coordinates {(1,0.15)};
      \addplot[blue]coordinates{(2, 0.2352)};
      %\addplot[blue]coordinates{(3,8.7e-06)};
      \addplot[blue]coordinates{(3, 0.05375 )};
      \addplot[blue]coordinates{(4,0.001024+0.001219 )};
      \addplot[blue]coordinates{(5,0.34556026)};
    \end{axis}
  \end{tikzpicture}
    %\caption{PCG} 
    \label{breakdown_pcg}
  \end{subfigure}
  \quad
  \begin{subfigure}[b]{0.45\linewidth}%\centering
    \begin{tikzpicture}[baseline = {(0,-0.88)}]
    \begin{axis}[
        xticklabel style={align=center},
        width=3in,height=2.52in,
        %height=2.8in,
        %scale only axis,
        %clip=false,
        %separate axis lines,
        %axis on top,
        %title={\large{Breakdown time at p = 5 for PCG}},
        %title style={yshift=0.75ex},
        ymajorgrids = true,
        xmin=0.5,
        xmax=4.5,
        xtick={1,2,3,4},
        x tick style={draw=none},
        xticklabels={\footnotesize v-cycle,\footnotesize CG \\[-3pt] \footnotesize matvec,\footnotesize {CG dot} \\[-3pt] \footnotesize product,\footnotesize PCG},
        ytick={0,0.1,0.2,0.3,0.4,0.5},
        ymin=0,
        ymax=0.5,
        ylabel={Time (s)},
        every axis plot/.append style={
          ybar=2pt,
          %x=1cm,
          bar width=0.5cm,
          bar shift=0pt,
          fill
        }
      ]
      \addplot[blue]coordinates {(1,0.34556026)};
      \addplot[blue]coordinates{(2, 0.03456)};
      \addplot[blue]coordinates{(3, 0.0017 )};
      \addplot[blue]coordinates{(4,0.41815885)};
    \end{axis}
  \end{tikzpicture}
    %\caption{v-cycle} 
    \label{breakdown_vcycle}
  \end{subfigure}
  \caption{Breakdown profiling at $p$ = 5 for the incompressible flow problem using GIAMG. (Left): profiling of v-cycle. ``First smooth" refers to the Chebyshev smoothing time at the first level; ``Smooth" refers to the time of Chebyshev smoothing at all levels; %"superLU" refers to the direct solving time at the coarsest level, 
  ``Residual" refers to the time of computing residual, which is projected to the next coarser level; ``PR transfer" refers to the time of residual restriction and error prolongation time; ``v-cycle" refers to the total time of one v-cycle. (Right): profiling of PCG. ``v-cycle" refers to the total time of one v-cycle; ``CG matvec" refers to the time of the additional matrix-vector multiplication other than v-cycle in the PCG algorithm; ``CG dot product" refers to the time of all additional dot product other than v-cycle in PCG; ``PCG" refers to total time of one PCG iteration.}
  \label{breakdown}
\end{figure}

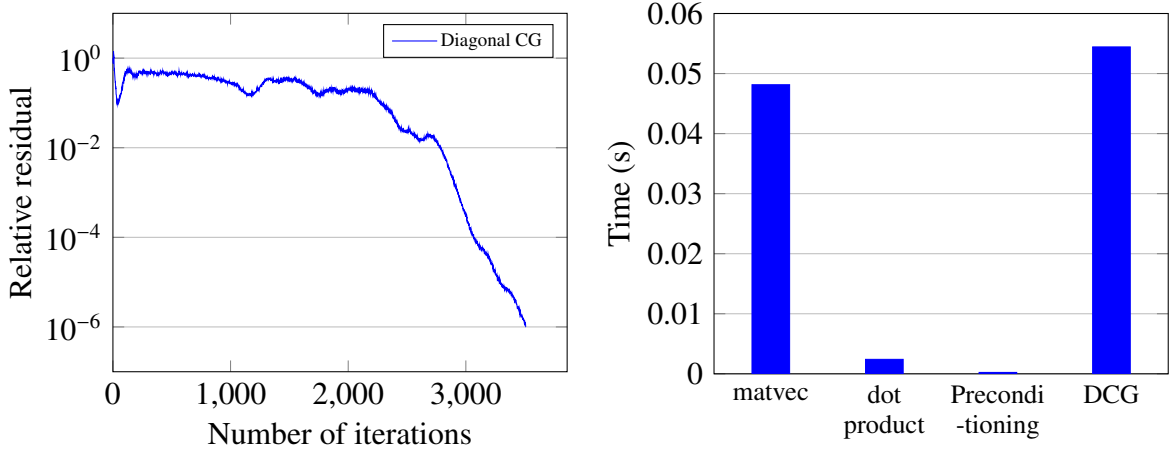
\begin{figure}[ht]
  \quad
  \begin{subfigure}[b]{0.45\linewidth}
    \begin{tikzpicture}%[scale=0.85]
      \begin{semilogyaxis}[width=3in,height=2.5in,
      ymajorgrids = true,
      xlabel={Number of iterations},ylabel={Relative residual},xmin=0,ymin=1e-7,ymax=10,legend entries={Diagonal CG},legend style={nodes={scale=0.6, transform shape}},legend pos=north east]
        \addplot [color=blue] table[] {nano/P5_dcg_paper.txt};
      \end{semilogyaxis}
    \end{tikzpicture}
    %\caption{Convergence} 
    \label{dcg_conv}
  \end{subfigure}
  \quad
  \begin{subfigure}[b]{0.45\linewidth}%\centering
    \begin{tikzpicture}[baseline = {(0,-1.08)}]
    \begin{axis}[
    yticklabel style={
        /pgf/number format/fixed,
        /pgf/number format/precision=5
        },
        scaled y ticks=false,
        xticklabel style={align=center},
        width=3in,
        height=2.51in,
        %scale only axis,
        ymajorgrids = true,
        %clip=false,
        %separate axis lines,
        %axis on top,
        %title={\large{DCG solve time at p = 5}},
        %title style={yshift=0.45ex},
        xmin=0.5,
        xmax=4.5,
        xtick={1,2,3,4},
        x tick style={draw=none},
        xticklabels={\footnotesize matvec, \footnotesize dot \\[-3pt] \footnotesize product, \footnotesize Precondi \\[-3pt] \footnotesize -tioning,\footnotesize DCG},
        ytick={0,0.01,0.02,0.03,0.04,0.05,0.06},
        ymin=0,
        ymax=0.06,
        ylabel={Time (s)},
        every axis plot/.append style={
          ybar=2pt,
          %x=1cm,
          bar width=0.5cm,
          bar shift=0pt,
          fill
        }
      ]
      \addplot[blue]coordinates{(1,4.813848e-02 )};
      \addplot[blue]coordinates{(2,2.406e-3)};
      \addplot[blue]coordinates{(3,2.286125e-04 )};
      \addplot[blue]coordinates{(4,5.443523e-02)};
    \end{axis}
  \end{tikzpicture}
    %\caption{Solve time} 
    \label{dcg_time}
  \end{subfigure}
  \caption{Breakdown profiling at $p$ = 5 for the incompressible flow problem using DCG solver. (Left): convergence. (Right): profiling. "matvec" refers to the time of matrix-vector multiplication in the CG algorithm; "dot product" refers to the time of all dot product in the CG algorithm; "Preconditioning" refers to the time of applying diagonal preconditioner within DCG; "DCG" refers to total time of one DCG iteration.}
  \label{dcg}
\end{figure}

\paragraph{Breakdown profiling}
We report breakdown profiling at $p = 5$ in Figure~\ref{breakdown}. We report the timings for different key parts of the v-cycle and the PCG solver. %We can see the Vcycle costs most of time in PCG solver. 
First, within the v-cycle, the pre- and postsmoothing (2 smoothing iterations) dominate the 
%, which is $70\%$, 
total v-cycle time, with the smoothing at the first (finest) level being the most expensive
%, which is $65\%$, 
of the total smoothing time. 
It is important to note that this will not be true if we do not employ $p$-coarsening, as the coarser levels will be denser and more expensive to evaluate.
Intergrid transfer (residual restriction and error prolongation) only costs a small amount of time. Note, the superLU direct solver costs a negligible amount of time (indicating a small enough coarsest level system), and therefore it is not shown in Figure~\ref{breakdown}. 
%(8.7e-06)
In addition, the v-cycle preconditioning dominates the total time of PCG solver, which is mostly due to the matrix-vector multiplication in the pre- and postsmoothing in v-cycle 
%(mostly at the finest level)
, and it makes the GIAMG much more expensive than a regular DCG solver for one PCG iteration. However, on the other hand, GIAMG has a much faster convergence. This motivates us to compare the total solve time of GIAMG with a regular DCG solver.        
%
% \begin{figure}
%     \centering
%     \includegraphics[width=0.45\linewidth]{figures/breakdown_P5_new1.eps}
%     %\vspace{-5mm}
%     \includegraphics[width=0.45\linewidth]{figures/breakdown_P5_new2.eps}
%     \caption{Breakdown profiling at p = 5. (Left): profiling of Vcycle. ``First smooth" refers to the Chebyshev smoothing time at the first level, ``Smooth" refers to the time of Chebyshev smoothing at all levels, ``superLU" refers to the direct solving time at the coarsest level, ``Residual" refers to the time of computing residual, which is projected to the next coarser level, ``PR transfer" refers to the time of residual restriction and error prolongation time, ``Vcycle" refers to the total time of one Vcycle. (Right): profiling of PCG. ``Vcycle" refers to the total time of one Vcycle, ``CG matvec" refers to the time of the additional matrix-vector multiplication other than Vcycle in the PCG algorithm, ``CG dot product" refers to the time of all additional dot product other than Vcycle in PCG, ``PCG" refers to total time of one PCG iteration.}
%     \label{breakdown}
% \end{figure}
%
We compare the GIAMG solver with the well-optimized DCG solver implemented in {\em Nektar++} for the solve time as shown in Figure~\ref{dcg}. First of all, from Figure~\ref{breakdown} (right) and Figure~\ref{dcg} (right), the time of matrix-vector multiplication looks similar in the two methods, which indicates it is a reasonable comparison between {\em Nektar++} and our GIAMG solver for the solve time even though they are two different implementations. Considering the dominant cost of v-cycle, our GIAMG is about $7$ times as great as DCG in terms of the solve time in one PCG iteration. However, since DCG solver is about ($3509$ iterations) $40$ times as great as our GIAMG solver in terms of the number of iterations for convergence, we are still about $6$ times as fast as the DCG solver in terms of the overall solve time.    
% \begin{figure}
%     \centering
%     \includegraphics[width=0.45\linewidth]{figures/comparison_dcg_p5_1.eps}
%     \includegraphics[width=0.45\linewidth]{figures/comparison_dcg_p5_2.eps}
%     \caption{Breakdown profiling at p = 5. (Left): convergence of DCG (Right): profiling of DCG. "Matvec" refers to the time of matrix-vector multiplication in the CG algorithm, "Preconditioning" refers to the time of applying diagonal preconditioner within DCG, "DCG" refers to total time of one DCG iteration.}
%     \label{dcg}
% \end{figure}

\begin{figure}[ht]
  \quad
  \begin{subfigure}[b]{0.45\linewidth}
    \begin{tikzpicture}%[scale=0.85]
      \begin{semilogyaxis}[width=3in,height=2.5in,
      ymajorgrids = true,
      xlabel={Number of iterations},ylabel={Relative residual},xmin=0,ymin=1e-7,ymax=10,legend entries={1 smooth iteration, 2 smooth iterations, 3 smooth iterations, 4 smooth iterations},legend style={nodes={scale=0.7, transform shape}},legend pos=north east]
        \addplot [very thick, smooth, color=red] table[] {nano/P5_new_1smooth_paper.txt};
        \addplot [very thick, smooth,color=blue] table[] {nano/P5_new_2smooth_paper.txt};
        \addplot [very thick, smooth,color=black] table[] {nano/P5_new_3smooth_paper.txt};
        \addplot [very thick, smooth,color=green] table[] {nano/P5_new_4smooth_paper.txt};
      \end{semilogyaxis}
    \end{tikzpicture}
    %\caption{Convergence} 
    \label{smooth_conv}
  \end{subfigure}
  \quad
  \begin{subfigure}[b]{0.45\linewidth}%\centering
    \begin{tikzpicture}[baseline = {(0,-1.05)}]
    \begin{axis}[
        yticklabel style={
        /pgf/number format/fixed,
        /pgf/number format/precision=5
        },
        x tick style={draw=none},
        xticklabel style={align=center},
        width=3in,height=2.53in,
        %height=2.75in,
        %scale only axis,
        %clip=false,
        %separate axis lines,
        %axis on top,
        %enlarge x limits=0.25,
        ymajorgrids = true,
        ylabel={Time (s)},
        %title={\large Time of varied smoothing at p = 5},
        %title style={yshift=0.2ex},
        % xmin=0.5,
        % xmax=4.5,
        % xtick={1,2,3,4},
        % x tick style={draw=none},
        % xticklabels={\small matvec, \small dot \\ product, \small Preconditioning,\small DCG},
        ytick={0,0.1,0.2,0.3,0.4,0.5,0.6,0.7},
        ymin=0,
        ymax=0.7,
        symbolic x coords={ \footnotesize {1 smooth} \\[-3pt] \footnotesize iteration, \footnotesize {2 smooth} \\[-3pt] \footnotesize iterations, \footnotesize {3 smooth} \\[-3pt] \footnotesize iterations, \footnotesize{4 smooth} \\[-3pt] \footnotesize iterations},
        xtick=data,
        x tick label style={
		/pgf/number format/1000 sep=},
	    enlarge x limits=0.12,
	    ybar=2pt,% configures `bar shift'
	    bar width=0.4cm,
	    legend pos=north west
      ]
      \addplot[fill=blue]
        coordinates {(\footnotesize {1 smooth} \\[-3pt] \footnotesize iteration,0.10911791) (\footnotesize {2 smooth} \\[-3pt] \footnotesize iterations,0.23525075)
         (\footnotesize{3 smooth} \\[-3pt] \footnotesize iterations,0.37654523) (\footnotesize {4 smooth} \\[-3pt] \footnotesize iterations,0.51612074)};
      \addplot[fill=red]
        coordinates {(\footnotesize {1 smooth} \\[-3pt] \footnotesize iteration,0.20676071) (\footnotesize {2 smooth} \\[-3pt] \footnotesize iterations,0.34556026) (\footnotesize {3 smooth} \\[-3pt] \footnotesize iterations,0.48383213) (\footnotesize {4 smooth} \\[-3pt] \footnotesize iterations,0.61968942)};
    %   \addplot[blue]coordinates{(2,2.406e-3)};
    %   \addplot[blue]coordinates{(3,2.286125e-04 )};
    %   \addplot[blue]coordinates{(4,5.443523e-02)};
    \legend{Smoothing, v-cycle}
    \end{axis}
  \end{tikzpicture}
    %\caption{Solve time} 
    \label{smooth_time}
  \end{subfigure}
  \caption{Comparison of different numbers of smoothing iterations at $p$ = 5 for the incompressible flow problem. Note, we apply the same number of iterations in both pre- and postsmoothing. (Left): convergence results. (Right) profiling of smoothing and v-cycle.}
  \label{smooth}
\end{figure}
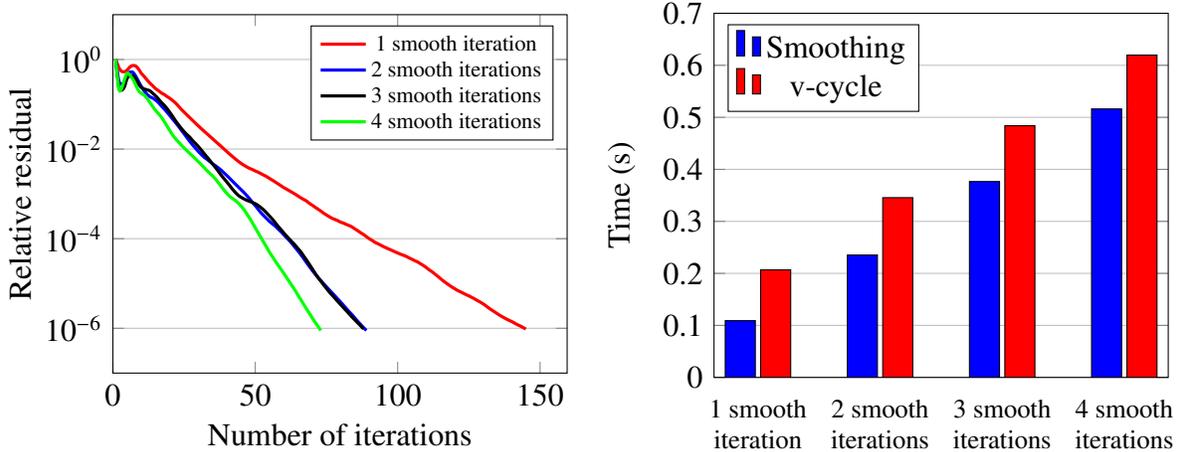

\paragraph{Number of Chebyshev smoothing iterations}
Since smoothing time dominates the overall time of v-cycle, it is the most important part that we want to control for achieving overall solving efficiency. Less smoothing iterations will reduce the smoothing time, but on the other hand, it may not be able to provide a reasonable smoothing effect for a good convergence. This trade-off motivates us to perform the smoothing experiments as shown in Figure~\ref{smooth}. We can see as increasing the number of smoothing iterations, the number of iterations for convergence generally decreases, while the time for smoothing as well as v-cycle increases. In addition, as we increase the number of smoothing iterations, the smoothing time becomes more and more dominant in the overall v-cycle time, further overshadowing other parts in the v-cycle. Specifically at $p = 5$, it can be concluded that the 2 smoothing iterations case is the optimal case in the experiments, and increasing number of smoothing iterations after it deteriorates overall solve performance, as the increased smoothing time overshadows the convergence gain. %The time costed in the rest of the PCG solver (i.e., other than Vcycle) is almost the same in all cases. 
Note, the effect of smoothing may be varied for different basis orders and different problems.      
% \begin{figure}
%     \centering
%     \includegraphics[width=0.45\linewidth]{figures/comparison_smooth_iteration_p5_new1.eps}
%     \includegraphics[width=0.45\linewidth]{figures/comparison_smooth_iteration_p5_new2.eps}
%     \caption{Comparison of different numbers of smoothing at p = 5. Note, we apply the same number iterations in both pre- and post-smoothing. (Left): convergence results. (Right) profiling of smoothing and vcycle. }
%     \label{smooth}
% \end{figure}
%

\begin{figure}[ht]
  \quad
  \begin{subfigure}[b]{0.45\linewidth}
    \begin{tikzpicture}%[scale=0.85]
      \begin{semilogyaxis}[width=3in,height=2.5in,
      ymajorgrids = true,
      xlabel={Number of iterations},ylabel={Relative residual},xmin=0,ymin=1e-7,ymax=10,legend entries={Coarsen by 2, Coarsen by 1},legend style={nodes={scale=0.8, transform shape}},legend pos=north east]
        \addplot [very thick, smooth,color=red] table[] {nano/P5_new_2smooth_diff_2_paper.txt};
        \addplot [very thick, smooth,color=blue] table[] {nano/P5_new_2smooth_diff_1_paper.txt};
      \end{semilogyaxis}
    \end{tikzpicture}
    %\caption{Convergence} 
    \label{p_diff_conv}
  \end{subfigure}
  \quad
  \begin{subfigure}[b]{0.45\linewidth}%\centering
    \begin{tikzpicture}[baseline = {(0,-1.0)}]
    \begin{axis}[
        yticklabel style={
        /pgf/number format/fixed,
        /pgf/number format/precision=5
        },
        xticklabel style={align=center},
        width=3in,height=2.54in,
        x tick style={draw=none},
        %height=2.78in,
        %scale only axis,
        %clip=false,
        %separate axis lines,
        %axis on top,
        %enlarge x limits=0.25,
        ymajorgrids = true,
        ylabel={ Time (s)},
        %title={\large Time of varied $p$-coarsening at p = 5},
        %title style={yshift=0.3ex},
        % xmin=0.5,
        % xmax=4.5,
        % xtick={1,2,3,4},
        % x tick style={draw=none},
        % xticklabels={\small matvec, \small dot \\ product, \small Preconditioning,\small DCG},
        ytick={0,0.1,0.2,0.3,0.4},
        ymin=0,
        ymax=0.4,
        symbolic x coords={\footnotesize{Coarsen by 1}, \footnotesize{Coarsen by 2}},
        xtick=data,
        x tick label style={
		/pgf/number format/1000 sep=},
	    enlarge x limits=0.4,
	    ybar=2pt,% configures `bar shift'
	    bar width=0.5cm,
	    legend style={nodes={scale=0.6, transform shape},legend columns=2},
	    legend pos=north east
      ]
      \addplot[fill=blue]
        coordinates {(\footnotesize{Coarsen by 1},0.14988872) (\footnotesize{Coarsen by 2},0.15050986)};
      \addplot[fill=red]
        coordinates {(\footnotesize{Coarsen by 1},0.23525075) (\footnotesize{Coarsen by 2},0.17325712)};
        \addplot[fill=yellow]
        coordinates {(\footnotesize{Coarsen by 1},0.053754023) (\footnotesize{Coarsen by 2},0.039391593)};
      \addplot[fill=green]
        coordinates {(\footnotesize{Coarsen by 1}, 0.34556026) (\footnotesize{Coarsen by 2},0.24670171)};
    %   \addplot[blue]coordinates{(2,2.406e-3)};
    %   \addplot[blue]coordinates{(3,2.286125e-04 )};
    %   \addplot[blue]coordinates{(4,5.443523e-02)};
    \legend{First smooth, Smooth, Residual, v-cycle}
    \end{axis}
  \end{tikzpicture}
    %\caption{Solve time} 
    \label{p_diff_time}
  \end{subfigure}
  \vspace{-0.45cm}
  \caption{Comparison of different $p$-coarsening strategies at $p$ = 5 for the incompressible flow problem. (Left): convergence results. (Right) profiling.}
  \label{p_diff}
\end{figure}
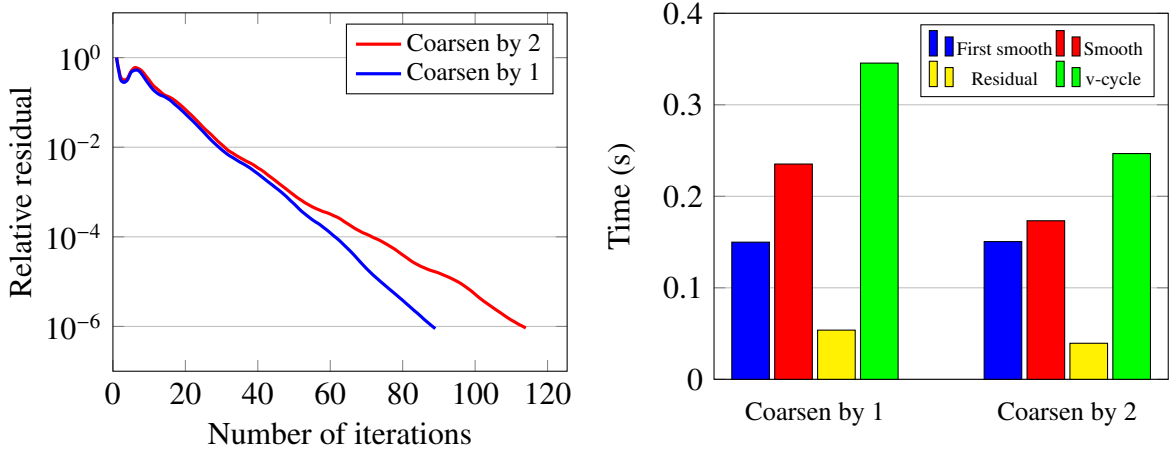

\begin{figure}
\centering
\begin{tikzpicture}
    \begin{axis}[width=3.5in,height=2.5in,
    xlabel={Number of nodes},ylabel={Time (s)},xmin=16,xmax=128,
    ymin=5e-4,ymax=5e-1,
    ymajorgrids = true,
    xmode=log,
    log basis x={2},
    ymode=log,
    log basis y={2},
    ymin=0.038,
    ymax=0.15,
    ytick={0.04,0.06,0.09,0.14},
    % yticklabel style={
    %     /pgf/number format/fixed,
    %     /pgf/number format/precision=5
    % },
    log ticks with fixed point,
    x tick label style={/pgf/number format/1000 sep=\,},
    legend entries={
    % Smooth,
    % %PR transfer,
    % Residual,
    v-cycle},legend style={nodes={scale=1, transform shape}},legend pos=north east]
        % \addplot [very thick, smooth, color=blue,mark=square*] coordinates {(16, 0.11239178) (32, 0.081753112) (64, 0.052832555) (128, 0.042145818) 
        % %(256, 0.044765702)
        % };
        % \addplot [very thick, smooth,color=red,mark=square*] coordinates {(16,0.00088117345+0.00079998759) (32,0.0006960639+0.000690573) (64,0.00061324494+0.00067188371) (128,0.00058021609+0.00069063298) %(256,0.00059581674+0.0007440699)
        %};
        % \addplot [very thick, smooth,color=green,mark=square*] coordinates {(16,0.015378507) (32,0.0090910094) (64,0.0059370155) (128,0.0054281645) 
        % %(256,0.0074454822)
        % };
        \addplot [very thick, smooth,color=blue,mark=square*] coordinates {(16,0.130786236044) (32,0.093022037315) (64,0.0605881666538) (128,0.0492063921197) %(256,0.0538109955012)
        };
    \end{axis}
\end{tikzpicture}
\caption{Strong scaling of v-cycle of GIAMG up to 128 nodes for the incompressible flow problem.}
    \label{fig:vcycle}
\end{figure}
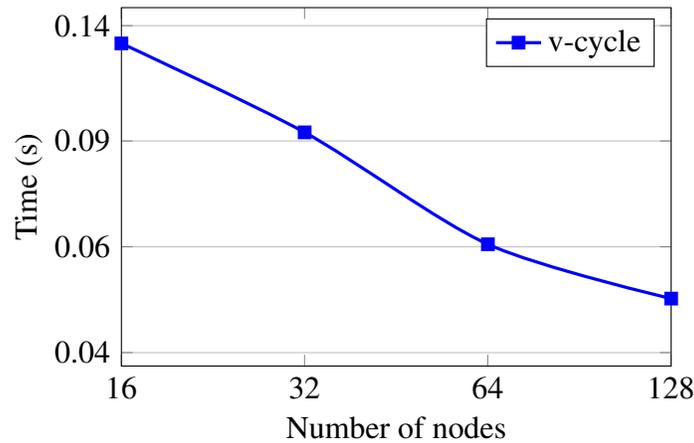

\paragraph{$p$-coarsening strategy}
Another key component that affects the GIAMG performance is the $p$-coarsening strategy. We test two coarsening strategies: (1), reducing basis order by 1 at each level (i.e., $5 \rightarrow 4 \rightarrow 3 \rightarrow 2 \rightarrow 1$); (2), reducing basis order by 2 at each level (i.e., $5 \rightarrow 3 \rightarrow 1$). The comparison is shown in Figure~\ref{p_diff}. The smoothing time is reduced when we coarsen the basis order by 2 at each level (the finest level smoothing time is still the same). The time saved here is essentially the second-finest level ($p$ = 4) and the fourth-finest level ($p$ = 2) compared with coarsening the basis order by 1. The residual computation time is also reduced as we have fewer levels to compute residuals when coarsening the basis order by 2. These result in a cheaper v-cycle for the more aggressive $p$-coarsening case. On the other hand, the number of iterations for convergence increases for the more aggressive $p$-coarsening, which is expected as it may introduce more aliasing error. However, the more aggressive $p$-coarsening case still achieves an overall better performance, indicating the time saved in the v-cycle overshadows the inferior convergence. 
%Again, the time costed in the rest of the PCG solver is almost the same in both cases. 
Note again, the effect of $p$-coarsening strategies may also be order and problem specific.
% \begin{figure}
%     \centering
%     \includegraphics[width=0.45\linewidth]{figures/compare_p_coasening_p5_new1.eps}
%     \includegraphics[width=0.45\linewidth]{figures/compare_p_coasening_p5_new2.eps}
%     \caption{Comparison of different coarsening strategies at p = 5. (Left): convergence result (Right): profiling of different parts in Vcycle.}
%     \label{p_diff}
% \end{figure}
%
\paragraph{Breakdown scalability}\label{scaling}
The strong scaling experiments are performed on the Frontera at the Texas Advanced Computing Center (TACC). The configuration of computing node used in the experiments is described below:
\begin{itemize}
 \item Node: \begin{itemize}
                    \item Processor: Intel Xeon Platinum $8280$ (``Cascade Lake").
                    \item Number of cores: $28$ per socket, $56$ per node.
                    \item Clock rate: $2.7Ghz$ (``Base Frequency").
                    \item ``Peak" node performance: $4.8TF$, double precision.
                   \end{itemize}
 \item Memory: $DDR-4$ memory, $192$GB/node.
 \item Network: Mellanox InfiniBand, HDR-$100$.
\end{itemize}
Figure~\ref{fig:vcycle} shows the strong scaling of v-cycle on $16-128$ nodes. It can be concluded that our GIAMG method is able to be scaled up to $128$ nodes ($7168$ MPI processes).
% \begin{figure}
%     \centering
%     \includegraphics[width=0.7\linewidth]{figures/strong_nano.pdf}
%     \caption{Strong Scaling up to 256 nodes}
%     \label{fig:strong_nano}
% \end{figure}

% \begin{figure}
%     \centering
%     \includegraphics[width=0.7\linewidth]{figures/strong_scaling_pCG001.png}
%     \caption{Strong Scaling up to 256 nodes}
%     \label{fig:pcg}
% \end{figure}

% \begin{figure}
%     \centering
%     \includegraphics[width=0.7\linewidth]{figures/strong_scaling_vcycle001.png}
%     \caption{Strong Scaling up to 256 nodes}
%     \label{fig:vcycle_bar}
% \end{figure}

% \begin{figure}
%     \centering
%     \includegraphics[width=0.7\linewidth]{figures/scaling.eps}
%     \caption{Strong scaling of v-cycle of GIAMG up to 256 nodes for the Navier--Stokes operator.}
%     \label{fig:vcycle}
% \end{figure}
%

%
% \begin{figure}
% \centering
% \begin{tikzpicture}
%     \begin{axis}[width=3.5in,height=2.5in,
%     %title={Strong scaling at p = 5},
%     xlabel={Number of nodes},ylabel={Time (s)},xmin=16,xmax=128,
%     ymin=7e-2,ymax=1,
%     ymajorgrids = true,
%     xmode=log,
%     log basis x={2},
%     ymode=log,
%     log basis y={2},
%     % yticklabel style={
%     %     /pgf/number format/fixed,
%     %     /pgf/number format/precision=5
%     % },
%     log ticks with fixed point,
%     x tick label style={/pgf/number format/1000 sep=\,},
%     legend entries={v-cycle},legend style={nodes={scale=0.7, transform shape}},legend pos=north east]
%         \addplot [very thick, smooth, color=blue,mark=square*] coordinates {(16, 0.65473416) (32, 0.34477519) (64, 0.20318642) (128, 0.11242437)};
%     \end{axis}
% \end{tikzpicture}
% \caption{Strong scaling of v-cycle of GIAMG up to $128$ nodes for the Navier--Stokes operator.}
%     \label{fig:vcycle1}
% \end{figure}

\paragraph{Comparison with other multigrid solvers}\label{petsc_comp}
We present the comparisons of convergence and total solve time between GIAMG and other AMG solvers run from PETSc. Note, for this more complex system, making close comparisons is not as easy as for the Helmholtz operator in \S~\ref{helmholtz_p8}. We try to make the overall comparisons as reasonable as possible according to our best knowledge of using these packages from PETSc. Some key parameters used in each AMG approach are shown in Table~\ref{tab:ns_p5}. 
\begin{table}
    \centering
    \begin{tabular}{|p{23mm}|p{25mm}|p{32mm}|p{18mm}|p{25mm}|}
    \hline
        & Smoother /\newline Iterations & Coarsening \newline scheme & Hierarchy \newline levels & Coarsest level \newline size / nonzeros\\
    \hline
    GAMG & Chebyshev / 2&Default SA & 6& 240 / 57600\\
    \hline
    ML & Chebyshev / 2& MIS& 6& 736 / 35164\\
    \hline
    BoomerAMG &Chebyshev / 2 &Falgout &7 & 602 / 51884\\
    \hline
    GIAMG & Chebyshev / 2 &$p$-coarsening+SA& 7&623 / 51072 \\
    \hline
    \end{tabular}
    \caption{Some key parameters used in each multigird approach at $p = 5$ for the incompressible flow problem.}
    \label{tab:ns_p5}
\end{table}

% The convergence of each approach is shown in Figure~\ref{ns_p5} (left). 
% \begin{figure}
%     \centering
%     \includegraphics[width=0.45\linewidth]{figures/comp_ns_conv.eps}
%     \includegraphics[width=0.45\linewidth]{figures/comp_ns_timing.eps}
%     \caption{Comparison between AMG methods at p = 5. (Left): convergence result (Right): total solve time}
%     \label{ns_p5}
% \end{figure}
%
\begin{figure}[ht]
  \quad
  \begin{subfigure}[b]{0.45\linewidth}
    \begin{tikzpicture}%[scale=0.85]
      \begin{semilogyaxis}[width=3in,height=2.5in,
      xlabel={Number of iterations},xmin=0,ymin=1e-7,ymax=10,
      ymajorgrids = true,
      ylabel={Relative residual},legend entries={BoomerAMG,ML,GAMG,GIAMG},legend style={nodes={scale=0.55, transform shape}},legend pos=north east]
        \addplot [very thick, smooth,color=green] table[] {nano/boomeramg_paper.txt};
        \addplot [very thick, smooth,color=blue] table[] {nano/ml_paper.txt};
        \addplot [very thick, smooth,color=red] table[] {nano/gamg_paper.txt};
        \addplot [very thick, smooth,color=black] table[] {nano/giamg_paper.txt};
      \end{semilogyaxis}
    \end{tikzpicture}
    %\caption{Convergence} 
    \label{ns_p5_converge}
  \end{subfigure}
  \quad
  \begin{subfigure}[b]{0.45\linewidth}%\centering
    \begin{tikzpicture}[baseline = {(0,-1.02)}]
    \begin{axis}[
        width=3in,height=2.53in,
        xticklabel style={align=center},
        %height=2.8in,
        %scale only axis,
        %clip=false,
        %separate axis lines,
        %axis on top,
        ymajorgrids = true,
        xmin=0,
        xmax=5,
        xtick={1,2,3,4},
        x tick style={draw=none},
        xticklabels={\footnotesize Boomer \\[-3pt] \footnotesize{-AMG},\footnotesize ML,\footnotesize GAMG,\footnotesize GIAMG},
        ytick={0,200,400,600,800,1000},
        ymin=0,
        ymax=1000,
        ylabel={Total solve time (s)},
        every axis plot/.append style={
          ybar=2pt,
          %x=1cm,
          bar width=0.5cm,
          bar shift=0pt,
          fill
        }
      ]
      \addplot[blue]coordinates {(1,4.5642e+02-2.5601e+01)};
      \addplot[blue]coordinates{(2,6.1378e+02-4.2066e+02)};
      \addplot[blue]coordinates{(3,2.0774e+03-1.1935e+03)};
      \addplot[blue]coordinates{(4,0.31622618*147)};
    \end{axis}
  \end{tikzpicture}
    %\caption{Solve time} 
    \label{ns_p5_timing}
  \end{subfigure}
  \caption{Comparison between multigrid methods at p = 5 for the incompressible flow problem. (Left): convergence result. (Right): total solve time}
  \label{ns_p5}
\end{figure}
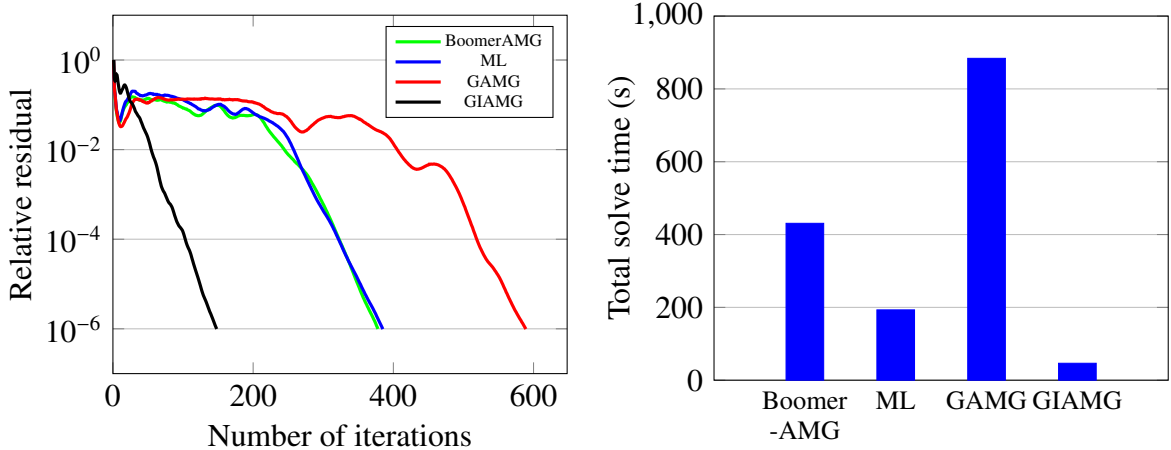
Our GIAMG method converges fastest among all the AMG packages under the relatively comparable conditions. In addition, it can be observed that other AMG solvers are stagnating at the relative residual around $0.1$ at the beginning for about 200 iterations, while GIAMG does not suffer this stagnation during convergence. This can also be observed when using DCG solver but with a much longer stagnation time as shown in Figure~\ref{dcg} (left). We speculate it is the $p$-coarsening in GIAMG that helps to achieve a significantly smoother convergence. After the stagnation period, all the AMG solvers seem to have a comparable convergence rate. The comparison of overall solve time is shown in Figure~\ref{ns_p5} (right). Again, our GIAMG solver performs the best among all the multigrid methods.    
% \begin{figure}
%     \centering
%     \includegraphics[width=0.7\linewidth]{figures/p5_timing.jpg}
%     \caption{Comparison of total solve time between multigrid methods at p = 5.}
%     \label{p5_timing_new}
% \end{figure}

\section{Conclusions and future work}~\label{conc}
In this paper, we proposed a geometrically informed algebraic multigrid (GIAMG) method, and developed an associated high-performance code for solving high-order finite element systems. %This method performs a p coarsening followed regular h coarsening. During p coarsening, additionall geometry information is required at the finest level. 
The GIAMG method sets up a grid hierarchy that includes $p$-coarsening at the top grids with assistance of minimal information of the geometry from outside.
We presented experiments for two 3D examples -- Helmholtz and incompressible flow problem. For the Helmheltz operator, we showed the convergence results from $p = 3$ up to $p = 10$. The convergence rate does not significantly change as the $p$ order is increasing. We also compared the convergence and solve time with some well-recognized AMG packages and a standard diagonally preconditioned PCG solver implemented in {\em Nektar++}. Our GIAMG performed the best among all the solvers. For the incompressible flow problem, we presented the convergence results from $p = 1$ up to $p = 5$. The convergence is also order-independent as we increase $p$ order. We further performed a detailed breakdown profiling of the GIAMG approach. The v-cycle dominates the PCG solving time, and pre- and postsmoothing take most time of the v-cycle. %With a similar matrix-vector multiplication efficiency, 
We are overall 6 times as fast as the diagonally preconditioned PCG solver. We further experimented on the effect of the number of smoothing iterations and $p$-coarsening strategies 
%(with fixed h coarsening) 
to the convergence and overall solve time. We also presented the strong scaling result of the GIAMG approach up to 128 nodes with 56 cores per node. Finally, we presented the comparisons of the GIAMG method with other AMG packages, and we again performed the best.  

Our future work is to identify and extract the geometric information required at the finest level by the GIAMG method purely from input matrices, making GIAMG to be a complete AMG approach. We are currently able to identify local elements, or cliques, from input matrices. However, we are still working on how to correctly determine the dofs which should be restricted to next coarser level within each element, and form the prolongation operator based on this coarsening.

\section*{Acknowledgements}~\label{acknowledgement}
The authors greatly appreciate the computing resources provided by the University of Utah Center for High-Performance Computing (CHPC) and the NSF-funded University of Texas Austin Texas Advanced Computing Center (TACC) machine (use of Frontera). We additionally thank Dr.~Andrew Davis and Dr.~Aleksander Dubas (UK Atomic Energy Authority) for supplying the geometry used in \S~\ref{NS}.
MR acknowledges support from  ARO W911NF-15-1-0222  (Program Manager Dr. Mike Coyle).
DM acknowledges support from the EPSRC Platform Grant PRISM under grant EP/R029423/1 and the ELEMENT project under grant EP/V001345/1. 
SX, RMK and HS acknowledge that this research was partially sponsored by ARL under Cooperative Agreement Number W911NF-12-2-0023. The views and conclusions contained in this document are those of the authors and should not be interpreted as representing the official policies, either expressed or implied, of ARL or the U.S. government. The U.S. government is authorized to reproduce and distribute reprints for government purposes notwithstanding any copyright notation herein. 

\section*{Appendix}~\label{App}
In this section, we describe how {\em Nektar++} orders elemental dofs. For a hexahedron element type, {\em Nektar++} defines the modal basis at each dof using tensor products of 1D basis functions. 1D basis functions $\phi$ are defined recursively using Jacobi polynomial as Equation~\ref{1d_basis}
\begin{equation}\label{1d_basis}
    \phi_{p}(a) =
    \begin{cases}
        \frac{1 - a}{2} & \text{if $p$ = 0} \\
        \frac{1 + a}{2} & \text{if $p$ = 1} \\
        \left(\frac{1-a}{2}\right)\left(\frac{1+a}{2}\right) P_{p-2}^{1,1}(a) & \text{otherwise}{,}
    \end{cases}
\end{equation}
where $p$ is the polynomial order, and $a$ is the 1D coordinate. $P_p^{\alpha,\beta}$ is the $p^{th}$ order Jacobi polynomial. For the linear case, we only have two modal basis functions in the element, $\frac{1-a}{2}$ and $\frac{1+a}{2}$, called the ``vertex modes", which coincides with the case of linear nodal basis. For high-order cases, higher order polynomials, the so-called ``edge modes" and ``interior modes", are introduced. The 3D basis $\Phi$ for an order p hexahedron element at each dof is defined as a tensor product of 1D basis functions
\begin{align}\label{modal_basis}
    \Phi_{(p+1)^2\times i+(p+1) \times j+k} = \phi_i \times \phi_j \times \phi_k
\end{align}
% in Algorithm~\ref{modal_basis}.
% %
% \begin{algorithm}[ht!]
%   \caption{Definition of modal basis}\label{modal_basis}
%   \begin{algorithmic}[1]
%   %\State {loop each element}
%     \Procedure{modal$\_$basis}{$i$,$j$,$k$}
%     %\State create $R$ and {\color{red}$n_R$\Comment{$n_R$ is the bin number}}
%         \For{$i=0:p$}
%             \For{$j=0:p$}
%                 \For{$k=0:p$}
%                     \State $\Phi_{(p+1)^2\times i+(p+1) \times j+k} = \phi_i \times \phi_j \times \phi_k$
%                 \EndFor
%             \EndFor
%         \EndFor
%     \EndProcedure
%   \end{algorithmic}
% \end{algorithm}
% %
Equation~\ref{modal_basis} also determines the ordering of the $(p+1)\times(p+1)\times(p+1)$ dofs within the hexahedron element in {\em Nektar++}, i.e., the $((p+1)^2\times i+(p+1)\times j+k)^{th}$ basis function is associated to the $((p+1)^2\times i+(p+1)\times j+k)^{th}$ dof in the elemental dof collection, implying that the order of this dof is $i$,$j$, and $k$ in the $x$, $y$ and $z$ directions of the coordinate system, respectively. For the spectral/$hp$ element basis definitions of other element types, we refer readers to~\cite{karniadakis2013spectral} for more detailed information.

% \section*{Appendix}~\label{App}

% The detailed choice of coarsening scheme and parameters used in Petsc options of each approach is shown as below. If not otherwise specified in the options, we use the default values for those parameters.
% Petsc GAMG options = "-ksp_type cg -pc_type gamg-pc_mg_galerkin external -pc_gamg_type agg -pc_gamg_agg_nsmooths 1 -mg_levels_ksp_type chebyshev -mg_levels_pc_type jacobi -mg_levels_ksp_max_it 3 -pc_gamg_threshold 0.015 -pc_gamg_sym_graph false -pc_gamg_square_graph 0 -ksp_monitor_true_residual -ksp_norm_type unpreconditioned -ksp_max_it 500 -ksp_rtol 1e-6 -ksp_converged_reason -ksp_view"; 

% ML option = "-ksp_type cg -pc_type ml -pc_mg_galerkin external -mg_levels_ksp_type chebyshev -mg_levels_pc_type jacobi -mg_levels_ksp_max_it 3 -pc_ml_maxNlevels 4 -pc_ml_Threshold 0.19 -pc_ml_CoarsenScheme Uncoupled -ksp_monitor_true_residual -ksp_norm_type unpreconditioned -ksp_max_it 500 -ksp_rtol 1e-6 -ksp_converged_reason -ksp_view"; 

% BoomerAMG option =  "-ksp_type cg -pc_type hypre -pc_hypre_type boomeramg -pc_hypre_boomeramg_max_levels 4 -pc_hypre_boomeramg_relax_type_all Chebyshev -pc_hypre_boomeramg_grid_sweeps_all 3 -pc_hypre_boomeramg_strong_threshold 0.28 -pc_hypre_boomeramg_coarsen_type HMIS -pc_hypre_boomeramg_agg_nl 2 -pc_hypre_boomeramg_agg_num_paths 3 -ksp_monitor_true_residual -ksp_norm_type unpreconditioned -ksp_max_it 500 -ksp_rtol 1e-6 -ksp_converged_reason -ksp_view -pc_hypre_boomeramg_print_statistics";

%%%%%%%%%%%%%%%%BIBLIOGRAPHY%%%%%%%%%%%%%%%%%
%%%%%%%%%%%%%%%%%%%%%%%%%%%%%%%%%%%%%%%%%%%%%
\bibliographystyle{unsrtnat}
\bibliography{sx}
%\bibliographystyle{plain}
%\bibliography{GraPTop}
% Graph of Diffuser Temperature Flow Rate Comparison

\end{document}